\documentclass[11pt]{amsart}
  \RequirePackage{amsmath}
\RequirePackage{amssymb} \RequirePackage{amsthm}
\RequirePackage{epsfig} 
\def\ep{\varepsilon}

\newcommand{\D}{\mathbb{D}}
\newcommand{\M}{\mathcal{M}}

\newcommand{\I}{\mathcal{I}}
\newcommand{\Fp}{\mathcal{F}_p}

\newcommand{\Fpf}{\mathcal{F}^{\phi}_p}
\newcommand{\Ftwof}{\mathcal{F}^{\phi}_2}
\newcommand{\Fqf}{\mathcal{F}^{\phi}_q}
\newcommand{\Fif}{\mathcal{F}^{\phi}_\infty}

\newcommand{\T}{\mathbb{T}}
\newcommand{\N}{\mathbb{N}}

\newcommand{\C}{\mathbb{C}}

\newcommand{\Z}{\mathbb{Z}}
\newcommand{\R}{\mathbb{R}}

\newcommand{\vp}{\vp}

\newcommand{\SSS}{\mathcal{S}}

\newcommand{\og}{\mathrm{O}}

\def\a{\alpha}       \def\b{\beta}        
            
     \def\om{\omega}      
              \def\f{\phi}
                  \def\z{\zeta}
                  \def\vf{\varphi}

\DeclareMathOperator{\supp}{supp}

\newtheorem{theorem}{Theorem}
\newtheorem{lemma}[theorem]{Lemma}
\newtheorem{proposition}[theorem]{\bf Proposition}

\newtheorem{corollary}[theorem]{\bf Corollary}

\newtheorem{lettertheorem}{Theorem}

\newtheorem{letterproposition}[lettertheorem]{Proposition}

\theoremstyle{definition}

\theoremstyle{remark}

\theoremstyle{remarks}

\newenvironment{pf}{\noindent{\emph{Proof. }}}{$\hfill\Box$ }
\newenvironment{Pf}{\noindent{\emph{Proof of}}}{$\hfill\Box$ }

\numberwithin{equation}{section}
\begin{document}
\title[Integral operators on weighted Fock spaces]
{Integral operators, embedding theorems and a Littlewood-Paley formula on weighted Fock spaces}

\author{Olivia Constantin}
\address{School of Mathematics, Statistics and Actuarial Science,
University of Kent,
Canterbury, Kent, CT2 7NF,
United Kingdom}
\email{O.A.Constantin@kent.ac.uk}
\address{
Faculty of Mathematics,
University of Vienna,
Nordbergstrasse 15, 1090 Vienna,
Austria}
\email{olivia.constantin@univie.ac.at}

\author{Jos\'e \'Angel Pel\'aez}

\address{Departamento de An´alisis Matem´atico, Universidad de M´alaga, Campus de
Teatinos, 29071 M´alaga, Spain} \email{japelaez@uma.es}

\thanks{The first author was supported in part by the FWF project
P 24986-N25. The second author was supported in part by the Ram\'on y Cajal program
of MICINN (Spain), Ministerio de Edu\-ca\-ci\'on y Ciencia, Spain,
(MTM2011-25502), from La Junta de Andaluc{\'i}a, (FQM210) and
(P09-FQM-4468).}
\date{\today}



\begin{abstract}
We obtain a complete characterization of the entire functions $g$ such that the integral operator $(T_ g f)(z)=\int_{0}^{z}f(\zeta)\,g'(\zeta)\,d\zeta$ is bounded or compact, on a large class of Fock spaces $\mathcal{F}^\phi_p$, induced by smooth radial
weights that decay faster than the classical Gaussian one. In some respects, these
 spaces turn out to be significantly different than the classical Fock spaces.
Descriptions of Schatten class integral operators are also provided.
\par En route, we prove a Littlewood-Paley formula for $||\cdot||_{\mathcal{F}^\phi_p}$ and we characterize the positive Borel measures for which $\mathcal{F}^\phi_p\subset L^q(\mu)$,\,$0<p,q<\infty$.
\par In addition, we also address the question of describing the subspaces of $\mathcal{F}^\phi_p$ that are invariant under the classical Volterra integral operator.

\end{abstract}
\maketitle


\section{Introduction}
\par \parindent 1em Let $\C$ be  the complex plane
and denote by $H(\C)$ the space of entire functions.
Given $\phi:[0,\infty)\rightarrow \R^+$  a twice continuously differentiable function and $0<p<\infty$,
we extend $\phi$
to $\C$ setting $\phi(z)=\phi(|z|)$, $z\in\C$.
We consider
the weighted Fock spaces,
$$\Fpf=\left\{f\in H(\C):\, ||f||^p_{\Fpf}=\int_{\C}|f(z)|^pe^{-p\phi(z)}\,dm(z) \right\}, $$
and
$$\Fif=\left\{f\in H(\C):\, ||f||_{\Fif}=\sup_{z\in\C}|f(z)|e^{-\phi(z)}\right\}, $$
where $dm$ denotes the Lebesgue measure in $\C$.
\par As usual we write $\Fp$ for the classical Fock spaces induced  by the standard function  $\f(z)=\frac{|z|^2}{2}$.
Moreover, for two real-valued functions $E_1,E_2$ we write $E_1\asymp
E_2$, or $E_1\lesssim E_2$, if there exists a positive constant $k$
independent of the argument such that $\frac{1}{k} E_1\leq E_2\leq k
E_1$, respectively $E_1\leq k E_2$.
\par We will deal with  the integral operator
$$T_g(f)(z)=\int_0^z f(\zeta)g'(\zeta)\,d\zeta$$
on  weighted Fock spaces $\Fpf$. This problem has been considered in \cite{ConsPAMS2012} for the classical Fock spaces proving that $T_g$ is bounded  only for
polynomials of degree $\le 2$. We expect that for a stronger decay of the weight $e^{-\phi}$ we shall find a wider class of symbols $g$ such that $T_g$ is bounded on $\Fpf$.

With this aim, we first  obtain a Littlewood-Paley type formula for $\mathcal{F}^\phi_p$ for a broad class of weights whose growth
ranges from logarithmic (e.g. $\phi(z)=a\log |z|,\, ap>2$) to highly exponential (e.g. $\phi(z)=e^{e^{|z|}}$).
Under some mild assumptions on $\phi$ (see Theorem \ref{Littlewood-Paley} below),
we prove that
\begin{equation}\label{L&P}
\int_\C |f(z)|^p e^{-p\phi(z)}\, dm(z)\asymp |f(0)|^p+
\int_\C |f'(z)|^p\left(\psi_{p,\phi}(z)\right)^p\,e^{-p\phi(z)}\,dm(z),
\end{equation}
 for any entire function $f$, where the {\it distortion function} $\psi_{p,\phi}$ satisfies
\begin{eqnarray}\label{compardist}
\psi_{p,\phi}(z)\asymp \frac{1}{\phi'(z)} \quad \hbox{ for } |z|\geq r_0,
\end{eqnarray}
for some $r_0>0$.
The above formula in the particular case of the classical Fock space  was
proven in \cite{{ConsPAMS2012}} using the explicit form of the reproducing kernel of $\mathcal{F}_2$.
The lack of precise information on the reproducing kernels for more general weights constrains us to employ a different method
based on estimates for integral means of entire functions and their derivatives.

\par
We then restrict our class of weights
 to rapidly increasing functions $\phi$. More precisely, we consider twice continuously differentiable functions $\phi:[0,\infty)\rightarrow \R^+$
 such that $\Delta\phi>0$ and
 \begin{eqnarray}\label{lap}
 \tau(z)\asymp\left\{
        \begin{array}{cl}
        1,  &  0\le |z|<1\\
 (\Delta \phi (|z|))^{-1/2}, &   |z|\ge 1
 \end{array}\right.,
 \end{eqnarray}
 where $\tau(z)$ is a radial positive differentiable function that
decreases to zero as $|z|\rightarrow\infty$ and $\lim_{r\rightarrow \infty} \tau'(r)=0$. Furthermore,
we  suppose that either there exists  a constant $C>0$
such that
$\tau(r)r^C$ increases for large $r$ or
 $$\lim_{r\to \infty}\tau'(r)\log\frac{1}{\tau(r)}=0.$$
 \par
 The class of  rapidly increasing functions $\phi$ will be denoted by $\I$.
It includes the power functions $\phi(r)=r^\alpha$ with $\alpha>2$ and exponential type functions such as
 $\phi(r)=e^{\beta r},$\, $\beta>0$ or $\phi(r)=e^{e^r}$.

 It turns out (see Theorem \ref{Littlewood-Paley} below) that for $\phi\in\I$, the Littlewood-Paley
 formula (\ref{L&P}) can be written in the form
 $$
 ||f||^p_{\Fpf} \asymp  |f(0)|^p+ \int_\C |f'(z)|^p\,\frac{e^{-p\f(z)}}{(1+\phi'(z))^p}\,dm(z).
 $$

Going further, we aim for a characterization of the Carleson measures for our Fock spaces which will  subsequently
be used to investigate the behavior of $T_g$.
This problem   was studied in \cite{carswell,ChoZhu2012,ortega} for the classical Fock space.
The Hilbert space case $p=2$ was considered by  Seip and Youssfi \cite{SeiYouJGA2011}  in several variables
for a wide class of radial weights.  However, our class of functions $\I$ does not completely overlap
with theirs, as shown by the examples $\phi(z)=|z|^m$ with $2<m<4$.
Let $D(a,r)$ be the Euclidean disc centered at $a$ with radius
$r>0$, and for simplicity, we shall write
$D(\delta \tau(a))$ for the disc $D(a,\delta
 \tau(a))$ with $\delta>0$ .
 \begin{theorem}\label{DCM}
Let  $\phi\in \mathcal{I}$ and let $\mu$ be a finite positive Borel
measure on $\C$.
\begin{enumerate}
\item[$(I)$] Let $0< p\leq q<\infty$.
\begin{enumerate}
 \item[$(a)$] The embedding
 $I_ d:\Fpf\rightarrow L^q(\mu)$ is bounded if and only if for
 some $\delta>0$ we have
\begin{eqnarray}\label{CMC}
K_{\mu,\phi}:=\sup_{a\in
\C}\frac{1}{\tau(a)^{2q/p}}\int_{D(\delta\tau(a))}\!\!e^{q\phi(z)}\,d\mu(z)<\infty.
\end{eqnarray}
Moreover, if any of the two equivalent conditions holds, then $$
K_{\mu,\phi} \asymp \|I_ d\|^q_{\Fpf\to L^q(\mu)}.$$

\item[$(b)$] The embedding
 $I_ d:\Fpf\rightarrow L^q(\mu)$ is compact if and only if for
 some $\delta>0$ we have
\begin{eqnarray}\label{CMC2}
\lim _{|a|\rightarrow
\infty}\frac{1}{\tau(a)^{2q/p}}\int_{D(\delta
\tau(a))}\!\!e^{q\phi(z)}\,d\mu(z)=0.
\end{eqnarray}
\end{enumerate}
\item[$(II)$] Let $0<q<p<\infty$. The following conditions are
equivalent:
\begin{enumerate}
\item[(a)] $I_ d: \Fpf\rightarrow L^{q}(\mu)$ is compact;
\item[(b)] $I_ d: \Fpf\rightarrow L^{q}(\mu)$ is bounded;

\item[(c)] For some $\delta>0$, the function
\begin{displaymath}
z\mapsto \frac{1}{\tau(z)^2} \int_{D(\delta\tau(z))}\!\!
e^{q\phi(\zeta)}\,d\mu(\zeta)
\end{displaymath}
belongs to $L^{\frac{p}{p-q}}(\C,dm)$.
\end{enumerate}
\end{enumerate}
\end{theorem}

This theorem has an interesting consequence. Recall that, for
the classical  Fock spaces $\Fp$, the following embeddings hold (see \cite{ZhuFock,tung})
$$\Fp\subseteq\mathcal{F}_q\quad \hbox{ for }0<p\le q<\infty.$$
For the particular choice $d\mu(z)=e^{-q\phi(z)}dm(z)$, Theorem \ref{DCM} shows under which conditions the Fock space $\Fpf$ is contained
in $\Fqf$  for $p\neq q$, where $p,q>0$.  
This never happens for rapidly increasing
functions $\phi\in \I$, as illustrated by the next corollary.

\begin{corollary}\label{co:nonested}
If $\phi\in\I$, the family of Fock spaces $\{\Fpf\}_{p>0}$ is not nested.  In fact, $\Fpf\setminus\Fqf\neq \emptyset$ and $\Fqf\setminus\Fpf\neq \emptyset$ for all $p,q>0$ with $p\neq q$.
\end{corollary}
\par This result shows a significant  difference between the weighted Fock spaces  $\Fpf$ with $\phi\in\I$, and the classical  Fock spaces $\Fp$,
which leads to additional technical difficulties in the study of $\Fpf$.

\par We then apply our Carleson embedding theorem together with the Littlewood-Paley formula to the study of the integral operator
$$T_g(f)(z)=\int_0^z f(\zeta)g'(\zeta)\,d\zeta.$$
The boundedness and compactness, as well as some
spectral properties (such as Schatten class membership) of $T_g$  acting on various spaces of analytic functions of the unit disc $\D$ in $\C$
have been extensively investigated (see e.g. \cite{cima,asis} for Hardy spaces, \cite{asi,PP,PelRat} for weighted Bergman spaces,
or the survey \cite{aa} and the references therein).

\begin{theorem}\label{th:Bq}
 Assume $g$ is an entire function, and $\phi\in\I$.
\begin{enumerate}
\item[$(I)$] For $0<p\le q<\infty$ we have
\begin{enumerate}
\item[$(a)$]
 $T_ g:\Fpf\rightarrow
\Fqf$ is bounded if and only if
\begin{eqnarray}\label{tgqbigp}
\sup_{z\in \C} \,\frac{|g'(z)|(\Delta\phi(z))^{\frac{q-p}{pq}}}{1+\phi'(z)}<\infty.
\end{eqnarray}
\item[$(b)$]  $T_ g:\Fpf\rightarrow
\Fqf$ is compact if and only if
\begin{displaymath}
\lim_{|z|\rightarrow \infty} \,\frac{|g'(z)|(\Delta\phi(z))^{\frac{q-p}{pq}}}{1+\phi'(z)}=0.
\end{displaymath}
\end{enumerate}
\item[$(II)$] Let $0<q<p<\infty$. The following conditions are
equivalent:
\begin{enumerate}
\item[$(a)$] $T_ g:\Fpf\rightarrow \Fqf$ is compact;

\item[$(b)$] $T_ g:\Fpf\rightarrow \Fqf$ is bounded;

\item[$(c)$] The function $\frac{g'(z)}{1+\phi'(z)}\in L^r(\C,dm)$, where $r=\frac{pq}{p-q}$.
\end{enumerate}
\end{enumerate}
\end{theorem}
\par Theorem \ref{th:Bq} shows a much richer structure of $T_g$ when acting on $\Fpf$ compared to the case of the classical Fock spaces,
where this operator is bounded only for polynomial symbols of degree $\le 2$ (see \cite{ConsPAMS2012}).
On the other hand,  Theorem \ref{DCM}, some other results in \cite{BDK,PP}, and the similar  techniques used to prove them,   show some analogies between Bergman spaces with rapidly decreasing weights
 $A^p_\om$ and $\Fpf$, $\phi\in\I$.
Indeed, the family of test functions considered to study the boundedness of $T_g$ on these spaces is basically the one introduced in \cite{BDK} to
characterize sampling and interpolation sequences (such problems for Fock spaces with nonradial weights $e^{-\phi}$ where $\Delta\phi$ is a doubling measure were considered in \cite{MarMasOrtGFA2003}).

 However, there are some fundamental differences as well:
 the weighted Bergman spaces are nested and only constant symbols induce bounded operators $T_g: A^p_\om\to A^q_\om$ when $0<p<q<\infty$ (see \cite{PP}),
 while Theorem \ref{th:Bq} illustrates a rich structure of $T_g:\Fpf\rightarrow \Fqf$
 for $0<p<q<\infty$ (see Section \ref{examples} for a further analysis of  Theorem \ref{th:Bq}).

Under a mild additional assumption on $\phi$,  we also provide a complete description of those
entire symbols $g$ for which the integral operator $T_g$
belongs to the Schatten $p$-class $\SSS_p(\Ftwof)$.
\begin{theorem}\label{th:Schatten}
Let $g$ be an entire function and assume $\phi\in \I$ satisfies
\begin{eqnarray}\label{CW5}
\sup_{r>r_0}\frac{-\tau''(r)\tau(r)}{\left(\tau(r)\phi'(r)\right)^2}<\infty,
\end{eqnarray}
for some $r_0>0$, where $\tau$ is as in (\ref{lap}).
\begin{enumerate}
\item[(a)] If $1<p<\infty$ then $T_ g\in \mathcal{S}_ p(\Ftwof)$
if and only if
$
 \,\frac{g'}{1+\phi'}\in L^p(\C,\Delta \phi\,dm).
$
\item[(b)] If $0<p\leq 1$ then $T_ g\in \mathcal{S}_ p(\Ftwof)$ if
and only if $g$ is constant.
\end{enumerate}
\end{theorem}
\par It is worth to comment that \eqref{CW5}  is a technical and non-restrictive condition, because for  natural examples from the class $\I$, the corresponding
function $\tau$ is convex. In addition, it is only used to prove the sufficiency for $p\ge 1$.

In Section \ref{invariant} we investigate the invariant subspaces of the Volterra operator $V:\Fpf\rightarrow\Fpf$, given by,
$$
V f(z)=\int_0^z f(\zeta) \, d\zeta.
$$
Notice that $V=T_g$ for $g(z)=z$.  For the weights $\phi(z)=|z|^m,\, m>2$ we obtain a complete characterization of the
invariant subspaces of $V$ showing that they are precisely the spaces
$$A^p_N=\overline{\hbox{Span}\{z^k\,:\, k\ge N \}}^{\Fpf}, \quad N\ge 0.$$
This illustrates the applicability of an abstract result stated below in Theorem \ref{invs} whose proof presents an interesting feature:
the characterization of the invariant subspaces of $V$ in the case $p\neq 2$ can be reduced to the case $p=2$ via Theorem \ref{th:Bq}.
The corresponding results for the classical Fock space can be found in \cite{{ConsPAMS2012}}. The fact that the weighted Fock
spaces $\Fpf$, $\phi\in\I$, are not nested requires a more involved approach  in comparison to the classical case.

We recall that a complete description of the invariant subspaces of $V$, when acting on various classical
spaces of analytic functions on the unit disc (Hardy spaces, standard
 Bergman spaces,
Dirichlet spaces) was obtained in \cite{ak}.
Furthermore, the real variable analogue has a long tradition that goes back to Gelfand and Agmon \cite{Ag,Gelfand}.

Finally, in Section \ref{bergman} we  point out that our approach on the Fock space also brings some improvements of the results in \cite{PP}
on Bergman spaces with rapidly decreasing weights. We first deduce a natural asymptotic estimate analogous to (\ref{compardist})
for the corresponding distortion function on the Bergman space. This observation leads us further to eliminate a hypothesis in
Theorem 2 from \cite{PP}, where a characterization of the boundedness and compactness of $T_g$ is provided. By doing this
we extend this characterization to a wider class of weights. In particular, we allow for a considerably faster decay, including
for example weights of the form
$$
\omega(z)=\exp({-e^{e^\frac{1}{1-|z|}}}).
$$

\par The paper is organized as follows. In Section \ref{spreliminaries} we present some preliminary results including the existence of a covering of $\C$ in terms of discs $D(\delta \tau(z))$. We  deal with the Littlewood-Paley formula \eqref{L&P} and some useful results concerning
the behavior of the function $\phi$  in Section \ref{sec:LP}. We prove
Theorem \ref{DCM}  in Section \ref{s3}, Theorem \ref{th:Bq} in Section \ref{s4} and Theorem \ref{th:Schatten} in Section \ref{Schatten}.
 Moreover, we  provide some  examples of functions in the class $\I$ in Section \ref{examples}. Section 8 is devoted to the study of invariant
 subspaces and in Section 9 we discuss Bergman spaces with rapidly decreasing weights.

\section{Preliminaries}\label{spreliminaries}
\subsection{Some technical tools}
In this section we present some facts, which are needed  to prove the main results, but which may also be of
independent interest.

\par Let $\tau$ be a positive function on $\C$. We say
that $\tau\in \mathcal{L}$ if there exist a constant  $c_
1>0$ such that
\begin{eqnarray}\label{EqB}
\big|\tau(z)-\tau(\zeta)\big |\leq c_ 1\,|z-\zeta|, \quad \textrm{
for } \quad z,\zeta\in \C.
\end{eqnarray}

\par Throughout this paper, we will always use the notation
$$m_{\tau}=\frac{\min\left(1,c^{-1}_1\right)}{4},$$
where  $c_1$ is the constant  in
\eqref{EqB}.

\begin{lemma}\label{tau}
Suppose that $\tau \in \mathcal{L}$,\, $0<\delta\le m_{\tau}$ and
$a\in\C$. Then,
$$\frac{1}{2}\,\tau(a)\le \tau(z)\le 2\,\tau(a)\quad\hbox{ if }\quad z\in D(a,\delta\tau(a)).$$
\end{lemma}
\begin{pf}
Note that, by condition \eqref{EqB} we have
\begin{displaymath}
\tau(a) \leq \tau(z)+c_ 1 |z-a|\leq \tau(z)+\frac{1}{4}\,\tau(a)
\quad \textrm{ if }\quad |z-a|\leq \delta \tau (a).
\end{displaymath}
Therefore $ \tau (a)\leq 2\,\tau(z)$  if $|z-a|\leq \delta \tau  (a)$.
 Similarly it can be proved that $\tau(z)\le 2\,\tau(a)$.
\end{pf}

We shall now sketch the proof of a covering lemma that is obtained by adapting an approach used
by Oleinik \cite{oleinik} for bounded domains to our setting.

\begin{lemma}\label{covering}
Assume $t: \C \to (0,\infty)$ is a continuous  function such that
\begin{eqnarray}\label{EqA}
|t(z)-t(\zeta)|\le \frac{1}{4} |z-\zeta|, \quad z,\zeta\in \C.
\end{eqnarray}
and $\lim_{|z|\rightarrow\infty} t(z)=0$.
Then there exists a sequence of points $\{z_j\}\in\C$ such that the following conditions are satisfied:

\begin{eqnarray*}
&(i)&
z_j\not\in D(z_k,t(z_k))\ for \  j\neq k;
\\
&(ii)&
\bigcup_{j\ge 1} D(z_j,  t(z_j))=\C;
\\
&(iii)&
\tilde D(z_j,  t(z_j))\subset D(z_j, 3t(z_j)), \hbox{ where }\tilde D(z_j, t(z_j))={\displaystyle \bigcup_{z\in D(z_j, t(z_j))}} D(z,t(z));
\\
&(iv)& \{D(z_j, 3t(z_j) )\}_{j \ge 1} \hbox{ is a covering of\, $\C$ of finite multiplicity.}
\end{eqnarray*}
\end{lemma}

\begin{pf}
We construct a sequence $\{z_j\}_{j \ge 1}$ inductively as follows: pick $z_1\in \C$ such that
$t(z_1)=\displaystyle\max_{z\in\C} t(z)$.  Provided $z_1,\,z_2,\, ...,\, z_{i-1}$ are chosen, we let $z_i$
be one of those the points in $C_i:=\C\setminus \Big(\bigcup_{k<i} D(z_k,t(z_k))\Big)$ such that $t(z_i)
=\displaystyle\max_{z\in C_i} t(z)$. This way, we obtain a sequence that satisfies condition $(i)$.
Conditions $(iii),(iv)$ are of local nature and therefore the proofs of $(3)-(4)$ from "Lemma of coverings"
in \cite[p. $233$]{oleinik} translate verbatim.

It remains to prove $(ii)$. We claim that $|z_n|\rightarrow\infty$ as $n\rightarrow\infty$. Indeed, if
this did not hold, $\{z_n\}_{n \ge 1}$ would possess a convergent subsequence $\{z_{n_k}\}_{k\ge 1}$. The
construction of $\{z_n\}_{n \ge 1}$ would then imply
$$
|z_{n_k}-z_{n_l}|\ge t(z_{n_k})\ge m>0, \quad k\neq l,
$$
since, by continuity, the positive function $t$ has a positive minimum on compacts.
This contradicts the convergence of $\{z_{n_k}\}_{k \ge 1}$ and the claim is proven.
Consequently, we have $t(z_n)\rightarrow 0$ as $n\rightarrow\infty$.
Now choose an arbitrary point $y\in\C$.  Since $t(y)>0$
there exists ${n_0}\ge 1$ with $t(y)>t(z_{n_0})$, which implies
$$
y\in \bigcup_{i<n_0} D(z_i,t(z_i)),
$$
and hence $(ii)$ holds.
\end{pf}
\par The next result can be obtained following the proof of \cite[Lemma 2.2]{PP}.

\begin{lemma}\label{le1}
Let $\f$ be a subharmonic  function
 and let
$\tau \in \mathcal{L}$ such that
 $\tau(z)^2 \,\Delta \f(z)\leq c_ 2$ for some constant $c_
 2>0$ and  $z\in \C\setminus\D$.  If $\beta\in \mathbb{R}$,  there exists a constant
 $M\geq 1$ such that
\begin{displaymath}
|f(a)|^p\,e^{-\beta\f(a)}\leq \frac{M}{\delta^2\,\tau(a)^2}\int_{D(
\delta \tau(a))}|f(z)|^p e^{-\beta\f(z)}\,dm(z),\quad a\in\C,
\end{displaymath}
for all $0<\delta\leq m_{\tau}$ and $f\in H(\C)$.
\end{lemma}

We note that if a function $\phi$ belongs to the class $\I$,
then its associated function $\tau(z)$ belongs to the class
$\mathcal{L}$. Thus Lemma \ref{le1} proves that for functions $\phi$ in
the class $\I$, the point evaluations $L_ a$ are bounded linear
functionals on $\Fpf$. Therefore, there are reproducing kernels
$K_ a\in \Ftwof$ with $\|L_ a\|=\|K_ a\|_{\Ftwof}$ and such that
\begin{displaymath}
L_ a f=\langle f, K_ a\rangle =\int_{\C} f(z)\,\overline{K_
a(z)}\,e^{-2\f(z)}\,dm(z), \quad f\in \Ftwof.
\end{displaymath}

Another consequence is that norm
convergence implies uniform convergence on compact subsets of
$\C$. It follows that also the space $\Fpf$ is
complete.\\
\subsection{Test Functions}
As a key tool for the study of the boundedness and compactness we use an  appropriate family of test functions constructed in \cite[Proposition $8.2$]{BDK}  (see also \cite{PP}).
\begin{letterproposition}\label{pr:BDK}
Let $\phi\in\I$ and  $R\ge 100$.  There exists $\eta(R)$ such that for every $a\in\C$ with $|a|\ge \eta(R)$, there exists an entire function $F_{a,R}$ such that
\begin{equation}\label{eq:testfunctionp1}
|F_{a,R}(\om)|e^{-\phi(\om)}\asymp e^{-\frac{|\om-a|^2}{4\tau^2(a)}}\asymp 1,\quad \om\in D(a,R\tau(a)),
\end{equation}

\begin{equation}\label{eq:testfunctionp2}
|F_{a,R}(\om)|e^{-\phi(\om)}\le C(\phi,R)\min\left\{1,\left[ \frac{\min\{\tau(a),\tau(\om)\}}{|\om-a|}\right]^{\frac{R^2}{2}}\right\},\quad \om\in \C.
\end{equation}
\end{letterproposition}

\begin{corollary}\label{co:testfunc}
Let $\phi\in \I$, $0<p<\infty$ and   $R> \max\{100,\frac{2}{\sqrt{p}}\}$, and $\eta(R)$ that from Proposition \ref{pr:BDK}. Then
\begin{enumerate}
\item[$(i)$] for $0< p<\infty$  the function $F_{a,R}$ in Proposition
\ref{pr:BDK} belongs to $\Fpf$ with $$\|F_{
a,R}\|_{\Fpf}^p\asymp \tau(a)^2, \quad \eta(R) \le |a|.$$

 \item[$(ii)$] the reproducing kernel $K_ a$ of
$\Ftwof$ satisfies the estimate
\begin{displaymath}
\|K_ a\|_{\Ftwof}^2\,e^{-2\phi(a)}\asymp \tau(a)^{-2}, \quad \eta(R)\le
|a|.
\end{displaymath}
\end{enumerate}
\end{corollary}
\begin{pf}
Let $a\in \C$ with $\eta(R)\le |a|$, and consider the functions
$F_{a,R}$. Write
$$R_k(a)=\left\{\om\in\C: 2^{k-1}R\tau(a)<|\om-a|\leq
2^kR\tau(a)\right\},\, \, k=1,2\dots.$$  Note that \eqref{eq:testfunctionp1}
gives
\begin{displaymath}
\int_{|\om-a|<R\tau(a)}|F_{a,R}(\om)|^pe^{-p\phi(\om)}\,dm(\om)\asymp \tau(a)^2,
\end{displaymath}
and, by \eqref{eq:testfunctionp2} and the fact that $R> \frac{2}{\sqrt{p}}$
\begin{eqnarray*}
\int_{|\om-a|>R\tau(a)}|F_{a,R}(\om)|^p\,e^{-p\phi(\om)}\,dm(\om)&\leq&
\sum_{k=1}^{\infty}\int_{R_k(a)}\!\!|F_{a,R}(\om)|^p\,e^{-p\phi(\om)}\,dm(\om)
\\
&\lesssim& \tau
(a)^{\frac{pR^2}{2}}\sum_{k=1}^{\infty}\int_{R_k(a)}\!\!\frac{dm(\om)}{|\om-a|^{\frac{pR^2}{2}}}
\\
&\lesssim &\sum_{k=1}^{\infty} 2^{-\frac{pR^2}{2} k} m\big(R_ k(a)\big )
\\
&\lesssim &\tau(a)^2.
\end{eqnarray*}
Therefore $F_{a,R}\in \Fpf$ with
$\|F_{a,R}\|_{\Fpf}^p\asymp \tau(a)^2$, which gives $(i)$.\\

\par  The use of Lemma \ref{le1} (with $\beta=2$) gives the
upper estimate of $(ii)$,
\begin{displaymath}
\|K_ a\|_{\Ftwof}^2\,e^{-2\phi(a)}\lesssim \tau(a)^{-2}.
\end{displaymath}
On the other hand, the functions $F_{a,R}$ obtained from the
previous proposition satisfy (by $(i)$) that $F_{a,R}\in \Ftwof$ with
$\|F_{a,R}\|_{\Ftwof}^2\asymp \tau(a)^2$, and by \eqref{eq:testfunctionp1} this
gives
\begin{displaymath}
|F_{a,R}(a)|^2\asymp e^{2\phi(a)}\asymp e^{2\phi(a)}\,\left(\tau(a)^2\right)^{-1}\,\|F_{a,R}\|_{\Ftwof}^2.
\end{displaymath}
Since $\|K_ a\|_{\Ftwof}=\|L_ a\|$, where $L_ a$ is the point
evaluation functional at the point $a$, this proves the lower
estimate of $(ii)$.
\end{pf}

\begin{proposition}\label{pr:fqp}
Let $\phi\in\I$, $0<p<\infty$ and
$R>\max{\left\{100, \frac{2}{\sqrt{p}},\,2\sqrt{p}\right\}}$. If $\eta(R)$ is the number
given in Proposition \ref{pr:BDK} and $\{z_ k\}\subset\C$ is the
sequence from
  Lemma \ref{covering}, the function
\begin{displaymath}
F(z)=\sum_{z_k: |z_k|\ge\eta(R)} \!\! a_ k \,\, \frac{F_{z_k, R}(z)}{\tau(z_ k)^{2/p}}
\end{displaymath}
belongs to $\Fpf$ for every sequence $\{a_ k\}\in \ell^p$.
Moreover,
\begin{displaymath}
\|F\|_{\Fpf}\lesssim \Big (\sum_ k|a_ k|^p \Big )^{1/p}.
\end{displaymath}
\end{proposition}
\begin{pf}
\par In what follows, we shall write
$$F(z)=\sum_{z_k: |z_k|\ge\eta(R)}  \!\! a_ k \,
\frac{F_{z_ k, R}(z)}{\tau(z_ k)^{2/p}}=\sum_{k}\, a_ k \,
\frac{F_{z_ k, R}(z)}{\tau(z_ k)^{2/p}},$$
and for simplicity we shall denote $\gamma=\gamma(R)=\frac{R^2}{4}$.
\par If $0<p\le 1$, then bearing in mind Corollary \ref{co:testfunc}, we
have that
\begin{equation*} \begin{split}\label{eq:pm1}
\|F\|^p_{\Fpf}&=\int_\C\left| \sum_ k a_ k \, \frac{F_{z_ k, R}(z)}{\tau(z_ k)^{2/p}}\right|^p\,e^{-p\phi(z)}\,dm(z)
\\ & \le
\sum_ k \frac{|a_ k|^p}{\tau(z_k)^2}\, \big \|F_{z_ k, R} \big
\|^p_{\Fpf}
\\ & \le C\sum_ k |a_ k|^p.
\end{split}\end{equation*}
\par If $p>1$, an application of H\"{o}lder's inequality
yields
\begin{equation}\label{eq:qp1}
|F(z)|^p \leq \sum_ k \frac{|a_ k|^p}{\tau(z_ k)^{2p}}\, |F_{z_ k, R}(z)|^{\frac{p(\gamma-p+1)}{\gamma}} \left (\sum_ k \tau(z_ k)^2\,
|F_{z_ k, R}(z)|^{p/\gamma}\!\right )^{p-1}\!\!\!\!.
\end{equation}
\par Now, we claim that
\begin{equation}\label{eq:qp2}
\sum_ k \tau(z_ k)^2\, |F_{z_ k, R}(z)|^{p/\gamma}\lesssim \tau(z)^2 e^{\frac{p\phi(z)}{\gamma}}\,.
\end{equation}
\par
In order to prove \eqref{eq:qp2}, fix $\delta_0\in (0,m_\tau)$ and observe that $t(z)=\delta_0\tau(z)$ satisfies the hypotheses of Lemma \ref{covering}. Using the
estimate \eqref{eq:testfunctionp1}, Lemma \ref{tau} and $(iv)$ of Lemma
\ref{covering}, we deduce that
\begin{eqnarray}
\begin{split}\label{eq:zk}
\sum_{\{z_k\in D(z,\delta_ 0 \tau(z))\}}\!\! \!\!\!\! \tau(z_
k)^2&\,|F_{z_ k, R}(z)|^{p/\gamma}
\\
&\lesssim e^{\frac{p\phi(z)}{\gamma}}\sum_{\{z_k\in D(z,\delta_ 0
\tau(z))\}}\!\!\!\!\!\! \tau(z_ k)^2
 \lesssim \tau(z)^2 e^{\frac{p\phi(z)}{\gamma}}\,.
\end{split}
\end{eqnarray}
\par On the other hand, an application of  \eqref{eq:testfunctionp2}
gives
\begin{displaymath}
\begin{split}
\sum_{\{z_k\notin D(z,\delta_ 0 \tau(z))\}}\!\! \!\!\!\! \tau(z_
k)^2 & \,|F_{z_ k, R}(z)|^{p/\gamma}
\\ & \lesssim
\tau(z)^{2p}e^{\frac{p\phi(z)}{\gamma}}\sum_{\{z_k\notin D(z,\delta_ 0
\tau(z))\}}\!\! \frac{\tau(z_ k)^2}{\,\,|z-z_k|^{2p}}
\\
&=\tau(z)^{2p}e^{\frac{p\phi(z)}{\gamma}}\,\sum_{j=0}^{\infty}\,\sum_{z_
k\in R_ j(z)} \frac{\tau(z_ k)^2}{\,\,|z-z_k|^{2p}},
\end{split}
\end{displaymath}
where
\begin{displaymath}
R_ j(z)=\left\{\zeta\in\C: 2^{j}\delta_ 0\tau(z)<|\zeta-z|\leq
2^{j+1}\delta_ 0\tau(z)\right\},\, \, j=0,1,2\dots
\end{displaymath}
\par By \eqref{EqB}, we deduce that,
for $j=0,1,2,\dots$, $$D(z_k,\delta_0\tau(z_k))\subset
D(z,5\delta_02^{j}\tau(z))\quad \textrm{ if }\quad z_k\in
D(z,2^{j+1}\delta_0\tau(z)).$$ This fact together with the finite
multiplicity of the covering (see Lemma \ref{covering}) gives
\begin{displaymath}
\sum_{z_ k\in R_ j(z)}\!\!\tau(z_ k)^2\lesssim
\,m\Big(D(z,5\delta_ 0 2^j\tau(z))\Big)\lesssim 2^{2j}\tau(z)^2.
\end{displaymath}
Therefore
\begin{equation*} \begin{split}\label{eq:ii}
\sum_{\{z_k\notin D(z,\delta_ 0 \tau(z))\}}\!\! \!\!\!\! \tau(z_
k)^2|F_{z_ k, R}(z)|^{p/\gamma} & \lesssim
\tau(z)^{2p}e^{\frac{p\phi(z)}{\gamma}}\sum_{j=0}^\infty\sum_{z_k\in R_
j(z)}\frac{\tau(z_ k)^2}{|z-z_k|^{2p}}
\\ & \lesssim e^{\frac{p\phi(z)}{\gamma}}\sum_{j=0}^\infty
2^{-2pj}\sum_{z_k\in R_ j(z)}\!\!\tau(z_ k)^2
\\ & \lesssim
\tau(z)^2 e^{\frac{p\phi(z)}{\gamma}}\sum_{j=0}^\infty 2^{(2-2p)j}
\\ & \lesssim \tau(z)^2 e^{\frac{p\phi(z)}{\gamma}},
\end{split}\end{equation*}
which together with \eqref{eq:zk}, proves  \eqref{eq:qp2}.\\
\par Now, joining \eqref{eq:qp1} and \eqref{eq:qp2}, we
obtain
\begin{equation*} \begin{split}
\|F\|^p_{\Fpf} &\leq \sum_ k \frac{|a_ k|^p}{\tau(z_
k)^{2p}}\int_{\C} |F_{z_ k, R}(z)|^{\frac{p(\gamma-p+1)}{\gamma}}\,\tau(z)^{2p-2}\,e^{-p\phi(z)\left(\frac{\gamma-p+1}{\gamma}\right)}\,dm(z).
\end{split}
\end{equation*}
So, it is enough to show that
\begin{eqnarray}\begin{split}\label{eq:qp3}
\int_{\C}|F_{z_ k, R}(z)|^{\frac{p(\gamma-p+1)}{\gamma}}\,\tau(z)^{2p-2}\,e^{-p\phi(z)\frac{\gamma-p+1}{\gamma}}\,\,dm(z)\lesssim
\, \tau(z_ k)^{2p}
\end{split}
\end{eqnarray}
to obtain the desired estimate
\begin{displaymath}
\|F\|^p_{\Fpf} \leq \sum_ k |a_ k|^p.
\end{displaymath}

\par It follows from \eqref{eq:testfunctionp1} and \eqref{EqB} that
\begin{eqnarray}\begin{split}\label{eq:qp4}
\int_{|z-z_ k|<\tau(z_ k)}&|F_{z_ k, R}(z)|^{\frac{p(\gamma-p+1)}{\gamma}}\,\tau(z)^{2p-2}\,e^{-p\phi(z)\frac{\gamma-p+1}{\gamma}}\,dm(z)
\\
& \asymp \int_{|z-z_ k|<\tau(z_ k)}\tau(z)^{2p-2}\,dm(z)\asymp
\tau(z_ k)^{2p}.
\end{split}
\end{eqnarray}
On the other hand, using \eqref{EqB}, it follows that
$$\tau(z)\leq C 2^j \tau(z_ k) \quad \textrm{ if } \quad |z-z_
k|<2^j \tau(z_ k).$$ Thus, since $\gamma>p$, bearing in mind
\eqref{eq:testfunctionp2}, we deduce that
\begin{displaymath}
\begin{split}
\int_{|z-z_ k|\geq\tau(z_ k)}\!\!&|F_{z_ k, R}(z)|^{\frac{p(\gamma-p+1)}{\gamma}}\,\tau(z)^{2p-2}\,e^{-p\phi(z)\frac{\gamma-p+1}{\gamma}}\,dm(z)
\\
& \lesssim \tau(z_ k)^{2p(\gamma-p+1)}\int_{|z-z_ k|\geq \tau(z_
k)}\frac{\tau(z)^{2p-2}}{|z-z_ k|^{2p(\gamma-p+1)}}\,dm(z)
\\
&\lesssim \tau(z_ k)^{2p(\gamma-p+1)}\!\sum_{j=0}^\infty
\int_{2^j\tau(z_k)\le |z-z_ k|<2^{j+1}\tau(z_
k)}\frac{\tau(z)^{2p-2}}{|z-z_ k|^{2p(\gamma-p+1)}}\,dm(z)
\\
&\lesssim \sum_{j=0}^\infty 2^{-2jp(\gamma-p+1)} \int_{2^j\tau(z_k)\le
|z-z_ k|<2^{j+1}\tau(z_ k)}\tau(z)^{2p-2}\,dm(z)
\\
&\lesssim \tau(z_ k)^{2p}\sum_{j=0}^\infty 2^{-2jp\,(\gamma-p)}
\lesssim \tau(z_ k)^{2p},
\end{split}
\end{displaymath}
which together with \eqref{eq:qp4} gives \eqref{eq:qp3}. This
finishes the proof.
\end{pf}

 \section{A Littlewood-Paley formula}\label{sec:LP}
Our aim will be to obtain a Littlewood-Paley formula for a large class of weighted Fock spaces $\Fpf$.
\par Let us introduce some useful notation. For an entire function $f$ and $0\le r<\infty$,
   we set
  \begin{equation*} \begin{split}
   M_p(r, f) &=\left (\frac{1}{2\pi }
  \int_{-\pi }^ \pi \vert f(re\sp {it})\vert^ p\, dt\right )^ {1/p},
  \quad 0<p<\infty
  ,\\
  I_p(r, f) & =M_p\sp p(r, f),\quad 0<p<\infty,
  \\   M_ \infty (r, f) & =\max_ {\vert z\vert =r}\vert f(z)\vert.
  \end{split}\end{equation*}

\par Inspired by considerations in \cite{Si} and \cite{PavP}, for  $\phi:[0,\infty)\rightarrow \R^+$   twice continuously differentiable function and $0<p<\infty$,
such that $\int_{0}^\infty s e^{-p\phi(s)}\,ds<\infty$,
 we
define
\begin{equation*} \label{eq:dis}
\psi_p(r)=\psi_{p,\f}(r)=\frac{\int_{r}^\infty s e^{-p\phi(s)}\,ds}{(1+r) e^{-p\phi(r)}}\quad 0\le r<\infty.
\end{equation*}  The function $\psi_{p,\f}$ will be called
the {\it{$p$-distortion function}}  of $\,\,\f$.
\par We consider the equivalence of the following conditions,
\begin{gather*}
\int_0^\infty M_q^p(r,f)re^{-p\phi(r)}\,dr\ <\infty,\\
\int_0^\infty M_q^p(r,f')(\psi_{p,\f}(r))^p\,re^{-p\phi(r)}\,dr\ <\infty,
\end{gather*}
for functions $f\in H(\C)$, where  $0<p<\infty$ and $0<q\le \infty$. The next theorem asserts that
this equivalence holds if $\phi$ fulfills the $K_p$-condition below, which is a rather weak assumption that is satisfied, for instance, whenever $\lim_{r\to\infty}r\phi'(r)=+\infty$.
\par\smallskip{\bf{Condition $\mathbf{K_p}$}}: The function $\f $ is differentiable and
there is a constant $K=K(p,\f)\in\R$ such that
\begin{eqnarray}\label{L}
\frac{\frac{d}{dr}\left(re^{-p\phi(r)}\right)\int_r^\infty se^{-p\phi(s)}\,ds}{r^2e^{-2p\phi(r)}} \le K,
\qquad 1\le r<\infty.
\end{eqnarray}

\begin{theorem}\label{th:LPf}
Let $0<p<\infty,$ $0< q\le \infty$. If $\,\f$ is a function satisfying
condition $K_p$, then
\begin{equation*}
\int_0^\infty M_q^p(r,f)re^{-p\phi(r)}\,dr \asymp |f(0)|^p+
\int_0^\infty M_q^p(r,f')(\psi_\f(r))^p\,re^{-p\phi(r)}\,dr
\end{equation*}
for any entire function $f$.
\end{theorem}
\par In particular, we obtain the following.
\begin{corollary}\label{co:LPf}
Assume that $0<p<\infty$ and $\,\f$ is a function satisfying the $K_p$-condition, then
\begin{equation*}
||f||^p_{\Fpf} \asymp |f(0)|^p+
\int_\C |f'(z)|^p\left(\psi_{p,\f}(z)\right)^p\,e^{-p\f(z)}\,dm(z)
\end{equation*}
for any entire function $f$.
\end{corollary}

\subsection{Proof of Theorem \ref{th:LPf}.}
We proceed in three steps: first we reformulate Theorem \ref{th:LPf} as Theorem \ref{th:LPf2}, then we present some preliminary material, and,
finally, we prove Theorem \ref{th:LPf2}.
\subsubsection{\bf{Reformulation}}
\par First, we note that for any fixed $R_0\in (0,\infty)$
\begin{eqnarray}
\label{1}
\int_0^\infty M_q^p(r,f)re^{-p\phi(r)}\,dr & \asymp&  \int_{R_0}^\infty M_q^p(r,f)re^{-p\phi(r)}\,dr \nonumber\\
&&
\\ \int_0^\infty M_q^p(r,f')(\psi_{p,\f}(r))^p\,re^{-p\phi(r)}\,dr & \asymp& \int_{R_0}^\infty M_q^p(r,f')(\tilde{\psi}_{p,\f}(r))^p\,re^{-p\phi(r)}\,dr\nonumber
\end{eqnarray}
where $\tilde{\psi}_{p,\f}(r)=\frac{\int_{r}^\infty s e^{-p\phi(s)}\,ds}{r e^{-p\phi(r)}}$ and the constants involved in (\ref{1}) depend on $\phi$, $p$ and $R_0$.
\par
Given a function
$\f,$ and $0<p<\infty,$ we define the function $\vf$ by
\begin{equation*}
\vf(r)^{-p}= p\int_r^\infty se^{-p\phi(s)}\,ds .
\end{equation*}
Since $\f$ is continuous
\begin{eqnarray}\label{eq:vfprime}
\vf(r)^{-p-1}\vf'(r)=  re^{-p\phi(r)},
\end{eqnarray}
and we  define the measure $dm_\vf$ on $[0,\infty)$ by
\begin{equation*} \label{mvf}
dm_\vf(r)=\frac{\vf'(r)}{\vf(r)}\,dr.
\end{equation*}
It is not difficult to see that condition $K_p$ is equivalent to
\begin{eqnarray}\label{M}
\sup_{1\le r<\infty}\frac{\vf''(r)\vf(r)}{\vf'(r)^2}\le M,
\end{eqnarray}
where $M\in\R$ is an appropriate constant.
\par Consequently, bearing in mind \eqref{1}, we can reformulate Theorem \ref{th:LPf} as follows.
\begin{theorem}\label{th:LPf2}
Let $0<p<\infty$ and $0< q\le \infty.$
 Let  $\vf:[1,\infty)\to \mathbb{R}\,$ be a differentiable, positive and increasing function
 such that $\lim_{r\to \infty}\vf(r)=\infty$. For each entire function $f$, we define
 $F_1(r):= \frac{M_q(r,f)}{\vf(r)}$, $F_2(r):=
\frac{M_q(r,f')}{\vf'(r)}$. If $\vf$ satisfies \eqref{M}, then
\begin{eqnarray}\label{eq:rf}
||F_1||^p_{L^p(dm_\vf)}\asymp |f(0)|^p+||F_2||^p_{L^p(dm_\vf)}
\end{eqnarray}
for each $f\in H(\C)$.
\end{theorem}
\subsubsection{\bf{Preliminary results}}

 We shall need
some lemmas.
\begin{lemma}\label{le1n}
If  $f\in H(\C),\,\,$ $0< q\le \infty,$ then there is a
constant $C_q$ such that
\begin{eqnarray}\label{conj}
M_q(r,f')\le C_q(\rho-r)^{-1}M_q(\rho,f),\quad
0\le r<\rho<\infty.
\end{eqnarray}
\end{lemma}
 \begin{pf}
\par Let $r$ and $\rho$  be as in the statement. Fix $R>\rho$ and define the analytic function in the unit disc $\D$,
$$f_R(z)=f(Rz),\quad z\in\D.$$
By \cite[Lemma $3.1$]{PavP}, there is $C_q>0$ such that
$$M_q(s,(f_R)')\le C_q\frac{M_q(t,f_R)}{t-s},\quad\text{for any $0<s<t<1$.}$$
Choosing $s=\frac{r}{R}$ and $t=\frac{\rho}{R}$, the proof is finished.
\end{pf}
\par The next  result can be proved following the lines of  \cite[Lemma $2$]{pac}.
\begin{lemma} \label{le2}
If  $f\in H(\C),\,\,$ $0<q\le \infty$ and $s=\min(q,1),$ then there is a constant $C_q$ such
that
\begin{equation*} \label{eq:le2}
M_q^s(\rho,f)-M_q^s(r,f)\le C_q(\rho-r)^sM_q^s(\rho,f'),
\qquad 0<r<\rho<\infty.
\end{equation*}
\end{lemma}
\par The following lemma can be obtained by standard techniques (see \cite[Lemma $3.4$]{PavP}), so its proof will be omitted.
\begin{lemma}\label{pac}
Let $\{A_n\}_{n=-1}^\infty$ be a sequence of complex numbers, $0<\gamma <\infty,$
$\alpha >0.\,$ Set
\begin{gather*}
Q_1=\sum_{n=-1}^\infty e^{-n\alpha}|A_n|^\gamma,\\
Q_2=|A_-1|^\gamma+\sum_{n=-1}^\infty
e^{-n\alpha}|A_{n+1}-A_n|^\gamma.
\end{gather*}
Then the  quantities $Q_1$ and $Q_2$ are equivalent in the sense that there is a
positive constant $C$
independent of $\{A_n\}_{n=-1}^\infty$ such that $(1/C)Q_1\le Q_2\le CQ_1$.
\end{lemma}
The next lemma, which will be a key tool for the proofs of
our main results, is essentially proved in \cite[Lemma $3.5$]{PavP}.
\begin{lemma}\label{2M}
 Let  $\vf:[1,\infty)\to \mathbb{R}\,$ be a differentiable, positive and increasing function
 such that $\lim_{r\to
\infty}\vf(r)=\infty$ and $\vf(1)=1.\,$   Define the sequence $\{r_n\}_{n=0}^\infty$
by
\begin{eqnarray}\label{r_n}
\vf(r_n)=e^n, \qquad n\ge 0.
\end{eqnarray}
If $\vf$ satisfies \eqref{M},\,  for every $n\ge 0,$
\[\begin{aligned}
\frac{\vf'(y)}{\vf'(x)}\le e^{2M},\quad r_n<x<y<r_{n+2}.
\end{aligned}
\]
\end{lemma}
\subsubsection{\bf{Proof of Theorem \ref{th:LPf2}.} }
\par Let $\{r_n\}_{n=0}^\infty$ be the sequence
defined by
\eqref{r_n}. We may assume without loss of generality that
$\vf(1)=1$.
\par
First, we assume  that $F_1\in L^p(dm_\vf)$. Then,
\begin{eqnarray}\label{dis}
\begin{aligned}
\infty&>\|F_1\|^p_{L^p(dm_\vf)} =\int_1^\infty M_q^p(r,f)\vf(r)^{-p-1}\vf'(r)\,dr\\
&\ge\sum_{n=0}^\infty M_q^p(r_n,f)\int_{r_n}^{r_{n+1}}
\vf(r)^{-p-1}\vf'(r)\,dr
\\ &=\sum_{n=0}^\infty M_q^p(r_n,f)\frac{\vf(r_{n})^{-p}-
\vf(r_{n+1})^{-p}}{p}
=C_p\sum_{n=0}^\infty M_q^p(r_n,f)e^{-np}.
\end{aligned}
\end{eqnarray}
On the other hand,
\begin{equation*}
\begin{aligned}
&\|F_2\|^p_{L^p(dm_\vf)}=
\int_1^\infty M_q^p(r,f')\vf'(r)^{1-p}\vf(r)^{-1}\,dr\\
&\le \sum_{n=0}^\infty M_q^p(r_{n+1},f')\int_{r_n}^{r_{n+1}}
\vf'(r)^{1-p}\vf(r)^{-1}\,dr
= \sum_{n=0}^\infty M_q^p(r_{n+1},f') \vf'(x_n)^{-p},
\end{aligned}
\end{equation*}
where $r_n<x_n<r_{n+1}$.
Here we have used the formula
$
\int_{r_n}^{r_{n+1}} \vf'(r)\vf(r)^{-1}\,dr =1.
$
Now, taking into account \eqref{conj} we obtain
\begin{eqnarray}\label{lab}
\|F_2\|^p_{L^p(dm_\vf)} \le C\sum_{n=0}^\infty
M_q^p(r_{n+2},f)(r_{n+2}-r_{n+1})^{-p} \vf'(x_n)^{-p}.
\end{eqnarray}
On the other hand, by Lagrange's theorem,
\begin{equation*}
\vf(r_{n+2})-\vf(r_{n+1})=(r_{n+2}-r_{n+1})\vf'(y_n),
\quad\text{where $r_{n+1}<y_n<r_{n+2}$},
\end{equation*}
whence
\begin{eqnarray}\label{LT}
r_{n+2}-r_{n+1} =(1-e^{-1})e^{n+2}\left(\vf'(y_n)\right)^{-1}.
\end{eqnarray}
Combining this with \eqref{lab} and Lemma \ref{2M}, we get
\begin{equation*} \begin{split}\label{lab1}
&|f(0)|^p+\|F_2\|^p_{L^p(dm_\vf)} \lesssim |f(0)|^p +\sum_{n=0}^\infty
M_q^p(r_{n+2},f)e^{-(n+2)p} \left(\frac{\vf'(y_n)}{
\vf'(x_n)}\right)^p
\\ & \lesssim |f(0)|^p+e^{2Mp}\sum_{n=0}^\infty
M_q^p(r_{n+2},f)e^{-(n+2)p}
 \lesssim \sum_{n=0}^\infty
M_q^p(r_{n},f)e^{-np},
\end{split}\end{equation*}
which together with \eqref{dis} gives the inequality
$|f(0)|^p+\|F_2\|^p_{L^p(dm_\vf)} \le C\|F_1\|^p_{L^p(dm_\vf)}$ in
\eqref{eq:rf}.
\par\medskip Now it will be proved the reverse inequality in \eqref{eq:rf}.
Assume  that $F_2\in
L^p(dm_\vf)$.
 We shall consider the case $q<1$, the proof for $q\ge 1$ is similar.
\par Let $q<1$ and $\gamma =p/q.$ Arguing as in  \eqref{dis} and choosing $r_{-1}=0$, we  get
\begin{equation*} \begin{aligned}
\|F_1\|^p_{L^p(dm_\vf)}&
\le C\sum_{n=-1}^\infty M_q^p(r_n,f)e^{-np}
=C \sum_{n=-1}^\infty A_n^\gamma e^{-np},
\end{aligned}
\end{equation*}
where $A_n= M_q^q(r_n,f)$.\,This together with Lemma \ref{pac} implies that
\[
\begin{aligned}
\|F_1\|^p_{L^p(dm_\vf)}
\le C|f(0)|^p +C \sum_{n=-1}^\infty
\big(M_q^q(r_{n+1},f)-M_q^q(r_n,f)\big)^{p/q} e^{-np}.
\end{aligned}\]
Hence, by Lemma \ref{le2}, we have that
\begin{eqnarray}\begin{split}\label{eq:fin}
\|F_1\|^p_{L^p(dm_\vf)}
\le C|f(0)|^p +C \sum_{n=-1}^\infty
M_q^p(r_{n+1},f')(r_{n+1}-r_n)^p e^{-np},
\end{split}\end{eqnarray}
Now, we use Lagrange's theorem, as in \eqref{LT}, to obtain
\begin{equation*} \label{F_1}
\begin{aligned}
\|F_1\|^p_{L^p(dm_\vf)}
&\lesssim |f(0)|^p + \sum_{n=-1}^\infty
M_q^p(r_{n+1},f')\vf'(x_n)^{-p},
\\ & \lesssim |f(0)|^p +|f'(1)|^p+\sum_{n=0}^\infty
M_q^p(r_{n+1},f')\vf'(x_n)^{-p}
\quad  \text{where $r_n<x_n<r_{n+1}$}.
\end{aligned}
\end{equation*}
On the other hand,
\[\begin{aligned}
\infty>&\|F_2\|^p_{L^p(dm_\vf)}=
\int_1^\infty M_q^p(r,f')\vf'(r)^{1-p}\vf(r)^{-1}\,dr\\
&\ge \sum_{n=0}^\infty M_q^p(r_{n+1},f')\int_{r_{n+1}}^{r_{n+2}}
\vf'(r)^{1-p}\vf(r)^{-1}\,dr\\
&= \sum_{n=0}^\infty M_q^p(r_{n+1},f') \vf'(y_n)^{-p},\quad
\text{where $r_{n+1}<y_n<r_{n+2}$}.
\end{aligned}
\]
Finally, using subharmonicity, \eqref{eq:fin},  Lemma \ref{2M} and  the above
inequality  we deduce
\begin{equation*} \begin{split}
\|F_1\|^p_{L^p(dm_\vf)}
&\le C|f(0)|^p +|f'(1)|^p+C \sum_{n=0}^\infty
M_q^p(r_{n+1},f')\vf'(x_n)^{-p}
\\ & \le
C|f(0)|^p +|f'(1)|^p+e^{2Mp}C\sum_{n=0}^\infty M_q^p(r_{n+1},f') \vf'(y_n)^{-p}
\\ & \le
C\left( |f(0)|^p +||F_2\|^p_{L^p(dm_\vf)}\right).
\end{split}\end{equation*}
The proof is complete.

\subsection{The distortion function}\label{distort}

We now consider a class of weights for which the statement of Corollary \ref{co:LPf} becomes more transparent.
\begin{lemma}\label{lhospital}
Assume $\phi:[0,\infty)\rightarrow \R^+$ is twice continuously differentiable and  there exists $r_0>0$ such that $\phi'(r)\neq 0$ for $r>r_0$.
Let $p>0$ and suppose
\begin{eqnarray}\label{eqn}
&&\lim_{r\rightarrow\infty} \frac{re^{-p\phi(r)}}{\phi'(r)}=0 \nonumber\\
&&\limsup_{r\rightarrow\infty} \frac{1}{r} \Bigl(\frac{r}{\phi'(r)}\Bigr)'<p\nonumber\\
&&\liminf_{r\rightarrow\infty} \frac{1}{r} \Bigl(\frac{r}{\phi'(r)}\Bigr)'>-\infty.
\end{eqnarray}
Then $\phi$ satisfies the condition $K_p$ and there exists $r_1>0$ such that
\begin{equation*}
\psi_p(r)\asymp \frac{1}{\phi'(r)} \hbox{ for } r\ge r_1,
\end{equation*}
where the involved constants might depend on $p>0$.
\end{lemma}

\begin{pf}
By hypothesis there is $\alpha<1$ and  $r_2\ge r_0$ such that
\begin{eqnarray}\label{r2}
\Bigl(\frac{r}{p\phi'(r)}\Bigr)'\le \alpha r\quad\text{on $[r_2,+\infty)$.}
\end{eqnarray}
 So
an integration by parts  on $(r_2,r]$ gives
\begin{equation*}\begin{split}
&\int_{r_2}^r se^{-p\phi(s)}\,ds= \int_{r_2}^r -p\phi'(s)e^{-p\phi(s)}\left( \frac{-s}{p\phi'(s)}\right)\,ds
\\& =\frac{-re^{-p\phi(r)}}{p\phi'(r)}+\frac{r_2e^{-p\phi(r_2)}}{p\phi'(r_2)}+\int_{r_2}^r \Bigl(\frac{s}{p\phi'(s)}\Bigr)' e^{-p\phi(s)}\,ds
\\ & \le \frac{-re^{-p\phi(r)}}{p\phi'(r)}+\frac{r_2e^{-p\phi(r_2)}}{p\phi'(r_2)}+\alpha\int_{r_2}^r  se^{-p\phi(s)}\,ds,
\end{split}\end{equation*}
that is
\begin{eqnarray*}
&\int_{r_2}^r se^{-p\phi(s)}\,ds\le\frac{1}{1-\alpha}\left(\frac{-re^{-p\phi(r)}}{p\phi'(r)}+\frac{r_2e^{-p\phi(r_2)}}{p\phi'(r_2)}\right),
\end{eqnarray*}
so taking the limit as $r\to \infty$ and bearing in mind \eqref{eqn}, we can assert that
 there is $C_1=C_1(p,\phi)>0$ such that for any $r\ge r_2$
\begin{eqnarray}\label{eq:ip2}
\int_r^\infty se^{-p\phi(s)}\le C_1\frac{re^{-p\phi(r)}}{\phi'(r)}.
\end{eqnarray}
In particular,  $\int_{0}^\infty se^{-p\phi(s)}\,ds<\infty$.
We note that $\phi$ satisfies the $K_p$ condition if
\begin{eqnarray}\label{eq:kp2}
\limsup_{r\to \infty}\left( \frac{\int_r^\infty se^{-p\phi(s)}\,ds}{r^2e^{-p\phi(r)}}- \frac{p\phi'(r)\int_r^\infty se^{-p\phi(s)}\,ds}{re^{-p\phi(r)}}\right)\le K<\infty.
\end{eqnarray}
Now, by (\ref{r2})
$$\frac{r}{\phi'(r)}=\int_{r_2}^r \Bigl(\frac{s}{\phi'(s)}\Bigr)'\,ds+\frac{r_2}{\phi'(r_2)}\le \frac{pr^2}{2}+\frac{r_2}{\phi'(r_2)},$$
which together with \eqref{eq:ip2} implies that
\begin{equation*}
\frac{\int_r^\infty se^{-p\phi(s)}\,ds}{r^2e^{-p\phi(r)}}\lesssim \frac{1}{r\phi'(r)} \lesssim \frac{p}{2}+\frac{r_2}{r^2\phi'(r_2)}<\infty,
\end{equation*}
so the first addend in \eqref{eq:kp2} is bounded. On the other hand, a straight-forward application of L'Hospital's rule gives
$$\liminf_{r\rightarrow\infty}\frac{p\int_r^\infty se^{-p\phi(s)}\,ds}{\frac{re^{-p\phi(r)}}{\phi'(r)}}\ge \liminf_{r\rightarrow\infty} \frac{1}{1-\frac{1}{pr} (\frac{r}{\phi'(r)})'}>-\infty,$$
consequently \eqref{eq:kp2} holds.
Finally, another  application of L'Hospital's rule gives
\begin{eqnarray*}
\liminf_{r\rightarrow\infty} \frac{1}{p-\frac{1}{r} (\frac{r}{\phi'(r)})'}\le
\liminf_{r\rightarrow\infty}\frac{\psi_p(r)}{(\phi'(r))^{-1}}\le \limsup_{r\rightarrow\infty}\frac{\psi_p(r)}{(\phi'(r))^{-1}} \le\limsup_{r\rightarrow\infty}\frac{1}{p-\frac{1}{r} (\frac{r}{\phi'(r)})'},
\end{eqnarray*}
which proves the lemma.
\end{pf}

\begin{lemma}\label{distortion}
Assume  $\phi:[0,\infty)\rightarrow \R^+$ is a twice continuously differentiable function such that $\Delta\phi>0\, ,(\Delta \phi (z))^{-1/2} \asymp \tau (z)$, where $\tau(z)$ is a radial positive function that
decreases to zero as $|z|\rightarrow\infty$ and $\lim_{r\rightarrow \infty} \tau'(r)=0$.
Then
\begin{eqnarray*}
&(a)& \lim_{r\rightarrow\infty} \frac{\phi'(r)}{r}=\infty.\\
&(b)& \lim_{r\rightarrow\infty}  \tau(r) \phi'(r)=\infty, \hbox{or, equivalently, } \lim_{r\rightarrow\infty} \frac{\phi''(r)}{(\phi'(r))^2}=0. \\
&(c)&
\hbox{For any $p>0$ we have }
\psi_p(r)\asymp \frac{1}{\phi'(r)+1}, \hbox{ for }
r\ge 0.\\
&(d)&  \liminf_{r\to \infty}\frac{r\Delta\phi(r)}{\phi'(r)}\ge C>0.
\end{eqnarray*}
\end{lemma}

\begin{pf}
A simple calculation shows that
\begin{eqnarray}\label{laplace}
\Delta \phi (r)= \phi''(r)+\frac{1}{r}\phi'(r)=\frac{(r\phi'(r))'}{r}.
\end{eqnarray}
By L'Hospital's rule we obtain
$$
\lim_{r\rightarrow\infty} \frac{r\phi'(r)}{r^2}=\lim_{r\rightarrow\infty}\frac{(r\phi')'}{2r}=\lim_{r\rightarrow\infty}\frac{\Delta \phi(r)}{2}=\infty,
$$
which proves $(a)$.

Let us now prove $(b)$. Taking into account (\ref{laplace}) we obtain
$$
\tau(r)\phi'(r)=\frac{\tau(r)}{r} \int_0^r s \Delta \phi (s)\, ds\gtrsim \frac{\tau(r)}{r} \int_0^r \frac{s}{\tau^2(s)}\, ds.
$$
Again by L'Hospital's rule we get
$$
\lim_{r\rightarrow\infty} \frac{r\tau^{-1}(r)}{\int_0^r \frac{s}{\tau^2(s)} \, ds}=
\lim_{r\rightarrow\infty} \frac {\tau(r)-r\tau'(r)}{r}=0,
$$
which implies $\lim_{r\rightarrow\infty}\tau(r)\phi'(r)=\infty$. By $(a)$ and relation (\ref{laplace}) this last
fact is equivalent to $\lim_{r\rightarrow\infty}\frac{ \phi''(r)}{(\phi'(r))^2}=0$.

Taking into account $(a)-(b)$ it is straight-forward to check that the hypotheses in Lemma \ref{lhospital} are satisfied, indeed
\begin{equation*}
\lim_{r\rightarrow\infty} \frac{1}{r} \Bigl(\frac{r}{\phi'(r)}\Bigr)'=0.
\end{equation*}
 and hence
$\psi_p(r)\asymp \frac{1}{\phi'(r)}$ for $r\ge r_0$. Since $\phi'\ge 0$ and $\lim_{r\rightarrow\infty } \phi'(r)=\infty$, we obtain
$(c)$.

We now turn to $(d)$. By \eqref{laplace}
$$\frac{r\Delta\phi(r)}{\phi'(r)}\asymp \frac{r\tau^{-2}(r)}{\phi'(r)}\asymp\frac{r^2\tau^{-2}(r)}{\int_0^r \frac{s}{\tau^2(s)} \, ds},$$
so  L'Hospital's rule implies
\begin{equation*}
\liminf_{r\to \infty}\frac{r^2\tau^{-2}(r)}{\int_0^r \frac{s}{\tau^2(s)} \, ds}\ge \liminf_{r\to \infty} 2\left(1-\frac{r\tau'(r)}{\tau(r)}\right)\ge 2,
\end{equation*}
and we are done.
\end{pf}

The previous considerations lead us to the following Littlewood-Paley formula.
\begin{theorem}\label{Littlewood-Paley}
Assume that $0<p<\infty$ and $\phi:[0,\infty)\rightarrow \R^+$ is twice continuously differentiable satisfying the hypotheses of Lemma \ref{lhospital}. Then we have
\begin{equation*}
||f||^p_{\Fpf} \asymp |f(0)|^p+ \int_\C |f'(z)|^p\left(\psi_{p,\phi}(z)\right)^p\,e^{-p\phi(z)}\,dm(z),
\end{equation*}
for any entire function $f$, where the {\it distortion function} $\psi_{p,\phi}$ satisfies
\begin{equation*}
\psi_{p,\phi}(z)\asymp \frac{1}{\phi'(z)} \quad \hbox{ for } |z|\geq r_0,
\end{equation*}
for some $r_0>0$. In particular, if $\phi\in\I$,   the following holds
\begin{equation*}
||f||^p_{\Fpf} \asymp  |f(0)|^p+ \int_\C |f'(z)|^p\,\frac{e^{-p\f(z)}}{(1+\phi'(z))^p}\,dm(z),
\end{equation*}
for any entire function $f$.
\end{theorem}

It is worth to  point out that the class of functions satisfying the hypotheses of Lemma \ref{lhospital}
is quite large, containing functions whose growth ranges from logarithmic (e.g. $\phi(r)=a\log (1+r),\, ap>2$) to highly exponential (e.g. $\phi(r)=e^{e^r}$).

The following result which is deduced from Lemma \ref{tau} shows that the distortion function is
"almost constant" on sufficiently small discs whose radii depend on $\Delta\phi$.
\begin{lemma}\label{almostct}
Assume that   $\phi:[0,\infty)\rightarrow \R^+$ is a twice continuously differentiable function such that
 $\Delta\phi>0\, ,(\Delta \phi (z))^{-1/2} \asymp \tau (z)$, where $\tau(z)$ is a radial positive function that
decreases to zero as $|z|\rightarrow\infty$ and $\lim_{r\rightarrow \infty} \tau'(r)=0$.
Then, there exists $r_0>0$ such that
\begin{eqnarray}\label{fi}
\phi'(a)\asymp\phi'(z), \quad z\in D(a, \delta \tau(a)),
\end{eqnarray}
for all $a\in\C$ with $|a|>r_0$. Moreover,
$$
1+\phi'(a)\asymp 1+\phi'(z), \quad z\in D(a, \delta \tau(a)),
$$
for all $a\in \C$.
\end{lemma}

\begin{pf}
Since $\phi$ is radial, it is enough to show that there exists $r_0>0$ such that for $a\ge r_0$
\begin{eqnarray}\label{ec}
\phi'(a)\asymp\phi'(r), \quad r\in (a-\delta \tau(a),a+\delta \tau(a)).
\end{eqnarray}
Recall that
$$
\phi'(r)=\frac{1}{r} \int_0^r s \Delta \phi (s)\, ds, \quad r>0.
$$
Our assumptions on $\tau$ imply
$$
a\asymp r, \quad r\in (a-\delta \tau(a),a+\delta \tau(a)),
$$
for $a>2\delta \max_{s\in[0,\infty)}\tau(s)$. Hence, proving (\ref{ec})
reduces to showing that
$$
\int_0^r s \Delta \phi (s)\, ds \asymp \int_0^a s \Delta \phi (s)\, ds, \quad r\in (a-\delta \tau(a),a+\delta \tau(a)).
$$
We have
\begin{eqnarray}\label{ec1}
\frac{\int_0^{a-\delta\tau(a)} \frac{s}{\tau^2 (s)}\, ds}{\int_0^a \frac{s}{\tau^2 (s)}\, ds}
\lesssim \frac{\int_0^r s \Delta \phi (s)\, ds}{\int_0^a s \Delta \phi (s)\, ds}
\lesssim \frac{\int_0^{a+\delta\tau(a)} \frac{s}{\tau^2 (s)}\, ds}{\int_0^a \frac{s}{\tau^2 (s)}\, ds}
\end{eqnarray}
By L'Hospital's rule and Lemma \ref{tau} we deduce
\begin{eqnarray*}
\limsup_{a\rightarrow\infty} \frac{\int_0^{a+\delta\tau(a)} \frac{s}{\tau^2 (s)}\, ds}{\int_0^a \frac{s}{\tau^2 (s)}\, ds}
&\lesssim&
\limsup_{a\rightarrow\infty} \frac{(1+\delta\tau'(a))(a+\delta \tau(a))\tau^{-2}(a+\delta \tau(a))}{a\tau^{-2}(a)}\\
&\sim& \limsup_{a\rightarrow\infty} \frac{(1+\delta\tau'(a))(a+\delta \tau(a))}{a}=1.
\end{eqnarray*}
Analogously one can show that
$$
\liminf_{a\rightarrow\infty} \frac{\int_0^{a-\delta\tau(a)} \frac{s}{\tau^2 (s)}\, ds}{\int_0^a \frac{s}{\tau^2 (s)}\, ds}\ge c>0.
$$
and (\ref{fi}) follows in view of (\ref{ec1}).
Since $\lim_{r\rightarrow\infty}\phi'(r)=\infty$ last assertion is now straight-forward.
\end{pf}

\section{Proof of Theorem \ref{DCM}}\label{s3}
 At the end of this section we provide a proof of Corollary \ref{co:nonested}.
 \par
Fix $R>\max\left\{100, \frac{2}{\sqrt{p}}, 2\sqrt{p}\right\}$, $\delta\in (0,m_{\tau})$ and
consider the covering $\{D(\delta \tau(z_j)\}$ given by Lemma \ref{covering} for $t(z)=\delta \tau(z)$.

\subsection {Proof of (I): boundedness}
 Suppose first
  that $I_ d:\Fpf\rightarrow L^q(\mu)$ is
bounded.
 For $a\in \C$ with $|a|\ge \eta(R)$, consider
the function $F_{a,R}$ obtained in Proposition
\ref{pr:BDK}.
By Corollary \ref{co:testfunc}, we have $\|F_{ a,R}\|_{\Fpf}^p\asymp
\tau(a)^2$. Then, using
(\ref{eq:testfunctionp1}) we get
\begin{displaymath}
\begin{split}
\int_{D(\delta\tau(a))}\!\! e^{q\phi(z)}\,d\mu(z)&\lesssim
\int_{D(\delta\tau(a))} |F_{ a,R}(z)|^q\,d\mu(z)\\&\le
 \int_{\C} |F_{a,R}(z)|^q\,d\mu(z)
\lesssim  \, \|I_ d\|^q_{\Fpf\to
L^q(\mu)}\,\tau(a)^{\frac{2q}{p}},
\end{split}
\end{displaymath}
which implies that $K_{\mu,\phi}\leq C \|I_ d\|^q_{A^p(w)\to L^q(\mu)}$.\\
\par Conversely, suppose that (\ref{CMC}) holds.
The idea of the proof goes back to \cite{oleinik}.
Bearing in mind  Lemma \ref{covering}, Lemma \ref{le1} and Lemma \ref{tau},
it follows that
\begin{eqnarray}\label{BEq1}
\begin{split}
\int_{\C} |f(z)|^q  d\mu(z)
&\leq
 \sum_ j \int_{D(\delta\tau(z_j))}|f(z)|^q\,d\mu(z)
\\
&\lesssim \sum_ j \int_{D(\delta\tau(z_
j))}\!\!\left(\frac{1}{\tau(z)^2}\int_{D(\delta\tau(z))}\!\!\!|f(\zeta)|^p e^{-p\phi(\zeta)}\,
dm(\zeta)\right)^{\frac{q}{p}} \!\!
e^{q\phi(z)}\,d\mu(z)
\\
&\lesssim \sum_ j \left (\int_{D(3\delta \tau(z_
j))}\!\!|f(\zeta)|^pe^{-p\phi(\zeta)}\,dm(\zeta)\right
)^{\frac{q}{p}}\int_{D(\delta\tau(z_
j))}\!\!\frac{e^{q\phi(z)}\,d\mu(z)}{\tau(z)^{\frac{2q}{p}}}\\
&\lesssim \, K_{\mu,\phi}\,\sum_ j \left (\int_{D(3\delta \tau(z_
j))}\!|f(\zeta)|^pe^{-p\phi(\zeta)}\,dm(\zeta)\right )^{\frac{q}{p}}
\end{split}
\end{eqnarray}
Now, using Minkowski inequality and the finite multiplicity $N$ of
the covering $\big\{D(3\delta\,\tau(z_ j))\big\}$ (see Lemma \ref{covering}), we have
\begin{displaymath}
\begin{split}
\int_{\C} |f(z)|^q \, d\mu(z)& \lesssim \,K_{\mu,\phi}\left (\sum_
j\int_{D(3\delta \tau(z_ j))} |f(\zeta)|^p\,e^{-p\phi(\zeta)}\,dm(\zeta)
\right )^{q/p}
\\
& \lesssim K_{\mu,\phi}\, N^{q/p}\,\|f\|_{\Fpf}^q,
\end{split}
\end{displaymath}
proving that $I_ d:\Fpf\rightarrow L^q(\mu)$ is continuous with
$\|I_ d\|^q_{\Fpf\to L^q(\mu)}\!\lesssim K_{\mu,\phi}$.

\subsection{ Proof of (I): compactness}

Suppose that (\ref{CMC2}) holds
and let $\{f_ n\}$ be a bounded sequence in $\Fpf$. By Lemma
\ref{le1}, Montel's theorem and  Fatou's lemma,
we may
extract a subsequence $\{f_{n_ k}\}$ converging uniformly on
compact sets of $\C$ to some function $f\in\Fpf$.
 Given $\ep>0$, fix
$r_0\in (0,\infty)$ with
\begin{eqnarray}\label{CEq2}
\sup_{|a|>r_
0}\frac{1}{\tau(a)^{2q/p}}\int_{D(\delta\tau(a))}e^{q\phi(z)}\,d\mu(z)<\ep.
\end{eqnarray}
Observe that there is $r'_ 0\in (0,\infty)$ with $r_0\le r'_0$ such that if a
point $z_k $ of the sequence $\{z_ j\}$ belongs to $\{|z|\le r_
0\}$, then $D(\delta\tau(z_k))\subset\{|z|\le r'_ 0\}$. So, take
$n_ k$ big enough such that $\sup_{|z|\leq r'_ 0}|f_{n_
k}(z)-f(z)|<\ep$. Then, setting $g_{n_ k}=f_{n_ k}-f$, and arguing
as in (\ref{BEq1}), it follows that
\begin{displaymath}
\begin{split}
\|g_{n_ k}\|_{L^q(\mu)}^q
&\leq \int_{|z|\leq r'_ 0} \!|g_{ n_
k}(z)|^q\,d\mu(z)+\!\sum_{|z_ j|>r_ 0}\int_{D(\delta\tau(z_ j))}
\!\!|g_{ n_ k}(z)|^q\,d\mu(z)
\\
&\le C \varepsilon +C\|g_{n_ k}\|_{\Fpf}^{q}\sup_{|z_ j|>r_
0}\frac{1}{\tau(z_ j)^{2q/p}}\int_{D(\delta\tau(z_
j))}\!\!e^{q\phi(z)}\,d\mu(z)
<  C\ep.
\end{split}
\end{displaymath}
This
proves that $I_ d:\Fpf\rightarrow L^q(\mu)$ is compact.

\par Conversely, suppose that $I_ d:\Fpf\rightarrow L^q(\mu)$ is
compact. Take
\begin{displaymath}
f_{a,R}(z)=\frac{F_{a,R}(z)}{\tau(a)^{2/p}}, \quad \eta(R)
\le |a|,
\end{displaymath}
where $\eta(R)$ and $F_{a,R}$ are obtained from Proposition \ref{pr:BDK}. By Corollary \ref{co:testfunc},
$$\sup_{|a|\ge \eta (R)} \|f_{a,R}\|_{\Fpf}\leq C<\infty$$ which together with the compactness of the
identity operator implies that $\{f_{a,R}: \eta (R) \le
|a|\}$ is a compact set in $L^q(\mu)$. Thus
\begin{eqnarray}\label{eq:c1}
\lim_{r\to \infty}\int_{r<|z|}|f_{a,R}(z)|^q\,d\mu(z)=0 \quad
\textrm {uniformly in }a.
\end{eqnarray}
\par On the other hand, if $\gamma=\frac{R^2}{4}$
  the estimate \eqref{eq:testfunctionp2} gives
\begin{displaymath}
|f_{a,R}(z)|^p \,e^{-p\phi(z)}\lesssim
\frac{\tau(a)^{2\gamma p-2}}{r^{2\gamma p}},\qquad |z|\leq r, \quad
|a|\geq 2r.
\end{displaymath}
Thus $f_{a,R}\rightarrow 0$ as $|a|\rightarrow \infty$ uniformly on
compact subsets of $\C$, which together with \eqref{eq:c1} implies
that $\lim_{|a|\rightarrow \infty}\|f_{a,R}\|_{L^q(\mu)}=0$.
Therefore, using the estimate \eqref{eq:testfunctionp1} of Proposition \ref{pr:BDK} we obtain
\begin{eqnarray*}
\tau(a)^{-2q/p}\int_{D(\delta\tau(a))}e^{q\phi(z)}d\mu(z)
&\lesssim&
\int_{D(\delta\tau(a))}|f_{a,R}(z)|^q\,d\mu(z)
\leq \|f_{a,R}\|_{L^q(\mu)}^q.
\end{eqnarray*}
Now let $|a|\rightarrow \infty$ above to complete the
proof.

\subsection{Proof of (II)}
\par The implication $(a)\Rightarrow (b)$ is obvious.
To prove that $(b)\Rightarrow (c)$, we use  an adaptation of
an argument due to Luecking (see \cite{L1}), where, instead of reproducing kernels, we employ the test functions  $F_{a,R}$.
For an arbitrary sequence $\{a_ k\}\in \ell^p$, consider the
function
\begin{displaymath}
G_ t(z)=\sum_{z_k: |z_k|\ge\eta(R)} a_ k \,r_ k(t)\, \frac{F_{z_ k,
R}(z)}{\tau(z_ k)^{2/p}}, \qquad 0<t<1,
\end{displaymath}
where $r_ k(t)$ is a sequence of Rademacher functions (see page
$336$ of \cite{L1}, or Appendix A of \cite{Duren1970}). By Proposition
\ref{pr:fqp}
and condition (b)
\begin{displaymath}
\int_{\C} |G_ t(z)|^q\,d\mu(z)\leq   C  \|G_t\|^q_{\Fpf} \leq C\left (\sum_ k |a_ k|^p\right
)^{q/p},\qquad 0<t<1.
\end{displaymath}
Integrating with respect to $t$ from $0$ to $1$, applying Fubini's
theorem, and invoking Khinchine's inequality (see \cite{L1}), we
obtain
\begin{eqnarray}\label{eq:K1}
\int_{\C} \left ( \sum_{z_k:\,|z_ k|\geq \eta(R)} \!\!|a_ k|^2
\frac{\big |F_{z_ k, R}(z)\big |^2}{\tau(z_ k)^{4/p}}\right
)^{q/2}\!\!\!\!\!d\mu(z)\leq C\!\left (\sum_ k |a_ k|^p\right
)^{q/p}\!\!\!\!\!\!\!.
\end{eqnarray}
\par
 If $\chi_ E(z)$ denotes the characteristic
function of a set $E$, bearing in mind the estimate \eqref{eq:testfunctionp1}, and
the finite multiplicity $N$ of the covering $\big\{D(3\delta\tau(z_
k))\big\}$ (see $(iv)$ of Lemma \ref{covering}), we have
\begin{displaymath}
\begin{split}
&\sum_{z_k:\,|z_ k|\geq \eta (R)}\frac{|a_k|^q}
{\tau(z_k)^{\frac{2q}{p}}}
\int_{D(3\delta\tau(z_k))}\!\!e^{q\phi(z)}\,d\mu(z)
\\
&\lesssim \sum_{z_k:\,|z_ k|\geq \eta(R)}\frac{|a_k|^q}
{\tau(z_k)^{\frac{2q}{p}}}
\int_{D(3\delta\tau(z_k))}|F_{z_k,R}(z)|^{q}\,d\mu(z)
\\
&=\int_{\C}\sum_{z_k:\,|z_ k|\geq \eta (R)}\frac{|a_k|^q}
{\tau(z_k)^{\frac{2q}{p}}}
|F_{z_k,R}(z)|^{q}\,\chi_{D\left(3\delta\tau(z_k)\right)}(z)\,d\mu(z)
\\&
\lesssim \max\{1,N^{1-q/2}\}\int_{\C}\left(\sum_{z_k:\,|z_ k|\geq \eta(R)} |a_ k|^2\frac{\big |F_{z_ k,R}(z)\big |^2}{\tau(z_
k)^{4/p}}\right)^{q/2}\!\!d\mu(z)
\end{split}
\end{displaymath}
This, together with \eqref{eq:K1} yields
$$\sum_{|z_ k|\geq \eta(R)}\frac{|a_ k|^q}
{\tau(z_ k)^{\frac{2q}{p}}} \int_{D(3\delta\tau(z_
k))}e^{q\phi(z)}\,d\mu(z)\le C \left (\sum_ k |a_ k|^p\right
)^{q/p}\!\!\!\!\!\!\!.$$
Then using the duality between $\ell^{ \frac{p}{q}}$ and
$\ell^{\frac{p}{p-q}}$ we conclude that
\begin{eqnarray}\label{qp31}
\sum_{ |z_ k|\geq \eta(R)}\left(\frac{1} {\tau(z_ k)^2}
\int_{D(3\delta\tau(z_ k))} \!\!
e^{q\phi(z)}\,d\mu(z)\right)^{\frac{p}{p-q}} \tau(z_ k)^2<\infty.
\end{eqnarray}
\par Note that there is $\rho_ 1\in (0,\infty)$, with $\eta(R)\le \rho_ 1$ such that if a
point $z_ k $ of the sequence $\{z_ j\}$ belongs to $\{|z|< \eta(R)\}$, then $D(\delta\tau(z_k))\subset\{|z|< \rho_ 1\}$. Therefore,
using Lemma \ref{tau}, $(ii)$ and $(iii)$ of Lemma
\ref{covering}, and \eqref{qp31} we deduce that
\begin{equation*}
\begin{split}\label{qp2}
&\int_{|z|\geq \rho_{1}}\!\!\left(\frac{1}
{\tau(z)^2}\int_{D(\delta\tau(z))} \!\!
e^{q\phi(\zeta)}\,d\mu(\zeta)\right)^{\frac{p}{p-q}} \! dm(z)
\\ & \le
\sum_{|z_ k|\geq \eta(R)}\int_{D(\delta\tau(z_k))}\left(\frac{1}
{\tau(z)^2}\int_{D(\delta\tau(z))} \!\!
e^{q\phi(\zeta)}\,d\mu(\zeta)\right)^{\frac{p}{p-q}} \! dm(z)
\\ & \lesssim
\sum_{|z_ k|\geq \eta(R)}\left(\frac{1} {\tau(z_
k)^2}\int_{D(3\delta\tau(z_ k))} \!\!
e^{q\phi(\zeta)}\,d\mu(\zeta)\right)^{\frac{p}{p-q}} \tau(z_
k)^2<\infty.
\end{split}
\end{equation*}
This, together with the fact that the integral
$$\int_{|z|<\rho_{1}}\!\!\left(\frac{1}
{\tau(z)^2}\int_{D(\delta\tau(z))} \!\!
e^{q\phi(\zeta)}\,d\mu(\zeta)\right)^{\frac{p}{p-q}} \! dm(z)$$ is
clearly finite, proves that $(c)$ holds.

\par Finally, we are
going to prove that $(c)$ implies $(a)$.
 It is
enough to prove that if $\{f_n\}$ is a bounded sequence in
$\Fpf$ that converges to $0$ uniformly on compact subsets of
$\C$ then $\lim_{n\to\infty}\|f_n\|_{L^q(d\mu)}=0.$ \par  Bearing in mind \eqref{EqA} for $t(z)=\delta\tau(z)$,
and the fact that $\tau(r)$ is a decreasing function with $\lim_{r\to\infty}\tau(r)=0$,
we assert that there is $r'_0$
\begin{eqnarray}\label{eq:z1}
D\left(\frac{\delta}{2}\,\tau(z)\right)\subset\left\{\zeta\in\C:\,
|\zeta|> \frac{r}{2}\right\},\quad\text{if}\quad |z|>r\ge r'_0.
\end{eqnarray}
\par On the other hand,
 it follows
from Lemma \ref{le1}  that
\begin{displaymath}
\begin{split}
|f_n(z)|^q &\leq C \frac{e^{q\phi(z)}}{\tau(z)^2}\int_{D
\left(\frac{\delta}{2}\tau(z)\right)}|f_n(\zeta)|^q\,
e^{-q\phi(\zeta)}\,dm(\zeta).
\end{split}
\end{displaymath}
Integrate with respect to $d\mu$, apply Fubini's theorem,  use
\eqref{eq:z1} and Lemma \ref{tau} to obtain
\begin{eqnarray}
\label{eq:z2}
&\displaystyle\int_{\left\{z\in\C:\, |z|> r\right\}}|f_n(z)|^q\,d\mu(z)\hspace{7cm}
\\
 & \le
C\!\displaystyle\int_{\left\{\zeta\in\C:\, |\zeta|>
\frac{r}{2}\right\}}|f_n(\zeta)|^q
e^{-q\phi(\zeta)}\left(\frac{1}{\tau(\zeta)^2}\displaystyle\int_{D(\delta\tau(\zeta))}\!\!e^{q\phi(z)}\,d\mu(z)\right)dm(\zeta),\,\,r\ge r'_0.\nonumber
\end{eqnarray}
By condition (c) for any fixed $\varepsilon>0$, there
is $r_ 0\ge r'_0\in (0,\infty)$, such that
$$\int_{\left\{\zeta\in\C:\, |\zeta|>
\frac{r_ 0}{2}\right\}}\left(\frac{1}
{\tau(\zeta)^2}\int_{D(\delta\tau(\zeta))}\!\!\!e^{q\phi(z)}d\mu(z)\right)^{\frac{p}{p-q}}\!\!\!dm(\zeta)
<\varepsilon^{\frac{p}{p-q}}.$$
 Then \eqref{eq:z2} and an application of
H\"{o}lder's inequality yields
\begin{eqnarray}
\label{eq:z3}
 &&\displaystyle\int_{\left\{z\in\C:\, |z|>r_ 0\right\}}|f_n(z)|^q\,d\mu(z) \\
 && \le C
\|f_n\|^{q}_{\Fpf}\left(\int_{\left\{\zeta\in\C:\, |\zeta|>
\frac{r_ 0}{2}\right\}}\left(\frac{1}
{\tau(\zeta)^2}\int_{D(\delta\tau(\zeta))}\!\!\!e^{q\phi(z)}d\mu(z)\right)^{\frac{p}{p-q}}\!\!\!dm(\zeta)
\right)^{\frac{p-q}{p}}\nonumber
\\ &&\le C\varepsilon.\nonumber
\end{eqnarray}
\par Moreover,
 we have
$\lim_{n\to\infty}\int_{ |z|\le r_ 0}|f_n(z)|^q\,d\mu(z)=0,$ which together with \eqref{eq:z3},
gives $\lim_{n\to\infty}\|f_n\|_{L^q(d\mu)}=0.$
 This completes the proof of Theorem \ref{DCM}.
 \medskip\\
  \begin{Pf} {\em{Corollary \ref{co:nonested}.}}
  Fix $p\in (0,\infty)$ and put  $\mu(z)=e^{-q\phi(z)}\,dm(z)$ in Theorem \ref{DCM}.
  For $q>p$, taking into account that $\tau$ decreases to $0$ as $a\to\infty$,
  we deduce that for any $\delta>0$
  $$\frac{1}{\tau(a)^{2q/p}}\int_{D(\delta\tau(a))}\!\!e^{q\phi(z)}\,d\mu(z)\asymp \tau(a)^{2(1-q/p)}\to\infty,\quad  \text{if $a\to\infty$},$$
  so by Theorem \ref{DCM} (I), $\Fpf  \not\subset \Fqf$.
   On the other hand, if $0<q<p$, using Theorem \ref{DCM} (II) and the fact that $f\equiv 1\notin L^{\frac{p}{p-q}}(\C,dm)$, we deduce that
  $\Fpf  \not\subset \Fqf$.
  \end{Pf}

\section{Proof of Theorem \ref{th:Bq}} \label{s4}
 Let us first notice that, by Theorem \ref{Littlewood-Paley},
 we have
\begin{eqnarray}\label{eq:z5}
\|T_ g f\|_{\Fqf}^q\asymp \int_{\C} |f(z)|^q \,|g'(z)|^q (1+\phi'(z))^{-q}\,e^{-q\phi(z)}\,dm(z),
\end{eqnarray}
for $q>0$ and for any entire function $f$.
This relation shows that the boundedness (compactness) of the
integration operator $T_ g:\Fpf\rightarrow \Fqf$ is equivalent
to the continuity (compactness) of the embedding $I_ d:\Fpf\rightarrow
L^q(\mu_{g,\phi})$ with
\begin{equation*}
d\mu_{g,\phi}(z)=|g'(z)|^q\,(1+\phi'(z))^{-q}\,e^{-q\phi(z)}\,dm(z).
\end{equation*}

\subsection{Proof of $\mathbf{(I)}$.}

Assume $0<p\leq q<\infty$ and let $\delta \in (0,m_{\tau})$.
By Theorem \ref{DCM}, the embedding $I_ d:\Fpf\rightarrow
L^q(\mu_{g,\phi})$ is continuous if and only if
\begin{eqnarray}\label{KE}
\sup_{a\in \C}\, \frac{1}{\tau(a)^{2q/p}}\int_{D(\delta\tau(a))}|g'(z)|^q\,(1+\phi'(z))^{-q}\,dm(z)<\infty.
\end{eqnarray}

Bearing in mind Lemma \ref{tau}, it is clear that \eqref{tgqbigp} implies \eqref{KE}. On the other hand, by Lemma \ref{almostct} and the subharmonicity of $|g'|^q$, we deduce
\begin{equation*} \begin{split}
|g'(a)|^q\,(1+\phi'(a))^{-q}\tau(a)^{2-2\frac{q}{p}} &
\lesssim (1+\phi'(a))^{-q}\tau(a)^{-2\frac{q}{p}}\int_{D(\delta\tau(a))}|g'(z)|^q\,\,dm(z)
\\ & \lesssim \tau(a)^{-2\frac{q}{p}}\int_{D(\delta\tau(a))}|g'(z)|^q(1+\phi'(z))^{-q}\,\,dm(z),
\end{split}\end{equation*}
and the proof of part $(a)$ of $(I)$ is complete.
\medskip\par
 Concerning the compactness part $(b)$ of $(I)$, note that, by part $(I)$ of
Theorem \ref{DCM},
the embedding $I_
d:\Fpf\rightarrow L^p(\mu_{g,\phi})$
is compact if and only if
\begin{displaymath}
\lim_{|a|\rightarrow
\infty }\frac{1}{\tau(a)^{\frac{2q}{p}}}\int_{D(\delta\tau(a))}|g'(z)|^q\,(1+\phi'(z))^{-q}\,dm(z)=0.
\end{displaymath}
Proceeding as in the boundedness part, we see that this is
equivalent to
\begin{displaymath}
\lim_{|a|\rightarrow \infty} |g'(a)|^q\,(1+\phi'(a))^{-q}\tau(a)^{2-2\frac{q}{p}}=0.
\end{displaymath}

\subsection{Proof of $\mathbf{(II)}$}
The equivalence $(a)\Leftrightarrow (b)$
follows from part $(II)$ of Theorem \ref{DCM}.
\\
\par Let us prove that $(b)\Rightarrow (c)$. From part $(II)$ of Theorem \ref{DCM} we deduce that
$(b)$ is equivalent to
\begin{displaymath}
\int_{\C} \left
(\frac{1}{\tau(z)^2}\int_{D(\tau(z))}\!\!\!\!\!|g'(\zeta)|^q (1+\phi'(\zeta))^{-q}\,dm(\zeta) \right
)^{\frac{p}{p-q}}\!\!\!\! dm(z)<\infty.
\end{displaymath}
Now,  by Lemma \ref{almostct} and the subharmonicity of $|g'|^q$,
\begin{displaymath}
\begin{split}
& \int_{\C} |g'(z)|^r (1+\phi'(z))^{-r}\,\,dm(z)
\\ \lesssim & \int_{\C} \left
(\frac{1}{\tau(z)^2(1+\phi'(z))^{q}}\int_{D(\tau(z))}\!\!\!\!\!|g'(\zeta)|^q \,dm(\zeta) \right
)^{\frac{p}{p-q}}\!\!\!\! dm(z)
\\ & \lesssim \int_{\C} \left
(\frac{1}{\tau(z)^2}\int_{D(\tau(z))}\!\!\!\!\!|g'(\zeta)|^q (1+\phi'(\zeta))^{-q}\,dm(\zeta) \right
)^{\frac{p}{p-q}}\!\!\!\! dm(z)<\infty.
\end{split}
\end{displaymath}
\par $(c)\Rightarrow (b)$. If $\frac{g'}{1+\phi'}\in L^r(\C,dm)$, then \eqref{eq:z5} and  H\"{o}lder's
inequality a  gives
\begin{displaymath}
\begin{split}
\|T_ g f\|_{\Fqf}^q &\lesssim \left (\int_{\C}
|f(z)|^p\,e^{-p\phi(z)}\,dm(z)\right )^{q/p}\left
(\int_{\C} |g'(z)|^r (1+\phi'(z))^{-r}\,\,dm(z)\right)^{q/r}
\\
&\lesssim \|\frac{g'}{1+\phi'}\|_{L^r(\C,\,dm)}^q\,\|f\|_{\Fpf}^q.
\end{split}
\end{displaymath}
Thus $T_ g:\Fpf\rightarrow \Fqf$ is bounded with $\|T_g
\|\lesssim \left\|\frac{g'}{1+\phi'}\right\|_{L^r(\C,dm)}$. This finishes the proof.\hfill $\Box$

\section{Schatten classes on $\Ftwof$}\label{Schatten}
\par Given a separable Hilbert space $H$, the
Schatten $p$-class of operators on $H$, $\mathcal{S}_
p(H)$, consists of those compact operators $T$ on $H$ with its
sequence of singular numbers $\lambda_ n$ belonging to $\ell^p$,
the $p$-summable sequence space.

If $\{e_ n\}$ is an orthonormal basis of a Hilbert space $H$ of
analytic functions in $\C$ with reproducing kernel $K_ z$, then
\begin{eqnarray}\label{RKformula}
K_ z(\zeta)=\sum_ n \langle K_z,e_n \rangle e_ n(\zeta)=\sum_ n e_ n(\zeta)\,\overline{e_ n(z)}
\end{eqnarray}
 Also, by
\eqref{RKformula} we have
\begin{displaymath}
\frac{\partial}{\partial\bar{z}}K_z(\zeta)=\sum_ n \overline{e'_
n(z)}e_ n(\zeta),\qquad z,\zeta\in \C.
\end{displaymath}
Thus Parseval's identity gives
\begin{eqnarray}\label{eqRKbasis}
\|K_ z\|_{H}^2=\sum_ n |e_ n(z)|^2 \quad\text{and}\quad \left|\left|\frac{\partial}{\partial\bar{z}}K_z\right|\right|^2_{H}=\sum_
n|e'_ n(z)|^2.
\end{eqnarray}

\par Now, we are going to give the proof of Theorem \ref{th:Schatten} on the description of the
Schatten classes $\mathcal{S}_ p:=\mathcal{S}_ p(\Ftwof)$. First
we consider the sufficiency part of the case $1< p<\infty$. For this we
need the following two lemmas.
The first one is a $L^\infty$ version of Theorem  \ref{th:LPf2} (see also  \cite[Theorem $2.1$]{PavCh}).

\begin{lemma}\label{th:LPf3}
Let $0<p\le\infty$ and  $\vf:[1,\infty)\to \mathbb{R}\,$ be a twice differentiable, positive and increasing function
 such that $\lim_{r\to \infty}\vf(r)=\infty$. If $\vf$ satisfies \eqref{M}, then for each entire function $f$, the following conditions are equivalent,
 \par (i)\, $M_p(r,f)=\og(\vf(r)),\quad r\to\infty$.
  \par (ii)\, $M_p(r,f')=\og(\vf'(r)),\quad r\to\infty$.
\end{lemma}
\begin{pf}
Let $\{r_n\}$ the sequence defined by $\vf(r_n)=e^n$.
\par $(i)\Rightarrow (ii)$.\,
By Lagrange's theorem for each $n\in \N$
\begin{eqnarray}\label{lagrange}
\vf(r_{n+1})-\vf(r_n)=\vf'(t_n)(r_{n+1}-r_n), \quad t_n\in (r_n,r_{n+1})
\end{eqnarray}
so by Lemma \ref{le1n} and Lemma \ref{2M}
\begin{equation*} \begin{split}
M_p(r_n,f') &\le  \frac{C_p}{r_{n+1}-r_n}M_p(r_{n+1},f)
  = \frac{C_p\vf'(t_n)}{\vf(r_{n+1})-\vf(r_n)}M_p(r_{n+1},f)
\\ & \le \frac{C_p\vf'(t_n)\vf(r_{n+1})}{\vf(r_{n+1})-\vf(r_n)}
 \le C \vf'(t_n)   \le C \vf'(r_n).
\end{split}\end{equation*}
Now, for $r\in[r_0,\infty)$  choose $n$ such that $r\in [r_n,r_{n+1})$. Then
$$ M_p(r,f')\le M_p(r_{n+1},f') \le C \vf'(r_{n+1})\le C \vf'(r),$$
where the last inequality follows from Lemma \ref{2M}.
\par $(ii)\Rightarrow (i)$.\, We assume that $0<p\le 1$ (the proof for $p>1$ is analogous). Bearing in mind Lemma \ref{le2}, \eqref{lagrange}
and Lemma \ref{2M}, we get
\begin{equation*} \begin{split}
M^p_p(r_{j+1},f)-M^p_p(r_{j},f) &\le C_p (r_{j+1}-r_j)^pM^p_p(r_{j+1},f')
 \le C (r_{j+1}-r_j)^p \vf'(r_{j+1})^p
\\ & \le C (r_{j+1}-r_j)^p\vf'(t_j)^p
 \le C e^{jp},
\end{split}\end{equation*}
so a summation  gives
$$M^p_p(r_{n+1},f)
\le C  \sum_{j=0}^n e^{jp} +M_p^p(r_
0,f)\le C e^{np}=C\vf(r_n)^p$$
and the proof follows.
\end{pf}

\begin{lemma}\label{le:sp1}
If $\phi\in\I$ satisfies \eqref{CW5}, then
$$\left|\left|\frac{\partial}{\partial\bar{z}}K_z\right|\right|_{\Ftwof}=\og
\left(\,\|K_ z\|_{\Ftwof}\phi'(|z|)\right), \quad |z|\to\infty.$$
\end{lemma}
\begin{pf}
\par Let $\{e_ n\}$ be the orthonormal basis of $\Ftwof$ given by
$$e_n(z)=z^n\delta^{-1}_n,\quad n\in\N,$$ where
$\delta^{2}_n=2\pi\int_0^\infty r^{2n+1}e^{-2\phi(r)}\,dr$. By Corollary \ref{co:testfunc}
 we have that
\begin{equation*}
\sum_{n=0}^\infty
r^{2n}\delta^{-2}_n=\sum_{n=0}^\infty|e_n(z)|^2=||K_z||^2_{\Ftwof}\asymp
\frac{e^{2\phi(r)}}{\tau^2(r)},\quad |z|=r.
\end{equation*}
 So, if we consider the entire function
defined by
$$f(z)=\sum_{n=0}^\infty z^n\delta^{-1}_n,$$ then
$M_2(r,f)=\left(\frac{1}{2\pi }\int_ {-\pi }\sp \pi \left \vert
f(re\sp {i\theta }) \right \vert \sp 2 d\theta
\right)^{\frac{1}{2}}\asymp \Phi(r),$ as $r\to \infty$, where
$$\Phi(r)=\frac{e^{\phi(r)}}{\tau(r)}.$$
\par By Lemma \ref{distortion}
\begin{equation*} \label{eq:sp4}
\Phi'(r)\asymp \Phi(r)\phi'(r), \quad\hbox{as }\ r\to \infty.
\end{equation*}
\par Moreover, in view of \eqref{CW5} and Lemma \ref{distortion}, a  calculation  shows that
\begin{equation*}
\limsup_{r\to \infty}\frac{\Phi''(r)\Phi(r)}{(\Phi'(r))^2}<\infty.
\end{equation*}
\par Thus by Lemma \ref{th:LPf3}
\begin{displaymath}
M_ 2(r,f')=\og \left (\Phi'(r)\right).
\end{displaymath}
Finally, since for $r=|z|$,
\begin{displaymath}
\begin{split}
\left|\left|\frac{\partial}{\partial\bar{z}}K_z\right|\right|^2_{\Ftwof}
= \sum_{n=0}^\infty|e'_n(z)|^2
 =\sum_{n=1}^\infty
n^2r^{2n-2}\delta^{-2}_n
 =M^2_2(r,f')
\end{split}
\end{displaymath}
we obtain
\begin{displaymath}
\left|\left|\frac{\partial}{\partial\bar{z}}K_z\right|\right|_{\Ftwof}=M_
2(r,f')=\og \left (\Phi'(r)\right)\asymp \Phi(r)\phi'(r)\asymp \|K_ z\|_{\Ftwof}\phi'(r),\quad r\to \infty.
\end{displaymath}
\end{pf}

\begin{proposition}\label{pr:Schsuf}
Let $g\in H(\C)$, $1<p<\infty$ and $\phi\in\I$ satisfying
\eqref{CW5}. If $\frac{g'}{1+\phi'}\in L^p(\C, \Delta\phi\,dm)$,  then
$T_g\in\mathcal{S}_p(\Ftwof)$.
\end{proposition}
\begin{pf}
\par
By
Theorem
\ref{Littlewood-Paley}, the inner product
$$\langle f,g\rangle_ *=f(0)\overline{g(0)}+\int_{\C}
f'(z)\overline{g'(z)}\,(1+\phi'(z))^{-2} e^{-2\phi(z)}\,dm(z)$$ gives a norm on
$\Ftwof$ equivalent to the usual one.
 If $1<p<\infty$, the operator $T_ g$
belongs to the Schatten $p$-class $\mathcal{S}_ p$ if and only if
\begin{displaymath}
\sum_ n \left|\langle T_g e_n,e_n\rangle_ *\right|^p<\infty
\end{displaymath}
for any orthonormal basis $\{e_ n\}$ (see \cite[Theorem $1.27$]{Zhu}).
Let $\{e_n\}$ be an orthonormal set of $\left(
\Ftwof,\langle,\rangle_*\right)$. Next, applying
Theorem \ref{Littlewood-Paley} for $p=1$ and $\breve{\phi}=2\phi$
we deduce
\begin{displaymath}
\begin{split}
1=  \|e_
n\|_{\Ftwof}^2= &\|e_
n^2\|_{\mathcal{F}^{\breve{\phi}}_1}\gtrsim
 \int_\C |e_n(z)e'_n(z)|(1+\breve{\phi}'(z))^{-1}e^{-\breve{\phi}(|z|)}\,dm(z)
\\ &\asymp  \int_\C |e_n(z)e'_n(z)|(1+\phi'(z))^{-1}e^{-2\phi(|z|)}\,dm(z)
\end{split}
\end{displaymath}
This together with
 H\"older's
 inequality yields
\begin{equation*} \label{eq:spbis}
\begin{split}
\sum_n\left|\langle T_g e_n,e_n\rangle_ *\right|^p  & \leq
\sum_n\left(\int_\C |g'(z)e_n(z)e'_n(z)|(1+\phi'(z))^{-2}e^{-2\phi(|z|)}\,dm(z)\right)^p
\\
& \lesssim \sum_ n \int_\C |g'(z)|^p|e_n(z)e'_n(z)|(1+\phi'(z))^{-(p+1)}e^{-2\phi(|z|)}\,dm(z)
\\
&=\int_\C |g'(z)|^p\left(\sum_n|e_n(z)e'_n(z)|\right)(1+\phi'(z))^{-(p+1)}e^{-2\phi(|z|)}\,dm(z),
\end{split}
\end{equation*}
and since $\|K_ z\|_{\Ftwof}^2 e^{-2\phi(|z|)}\asymp \Delta \phi(z)$ as $|z|\to\infty$ (see
Corollary \ref{co:testfunc}), the result will be proved if we are able to
show that
\begin{eqnarray}\label{eqsp:kl}
\sum_n|e_n(z)\,e'_n(z)|\lesssim \,\,\|K_
z\|_{\Ftwof}^2(1+\phi'(z)).
\end{eqnarray}
 To prove \eqref{eqsp:kl}, we use the Cauchy-Schwarz inequality
 to obtain
\begin{displaymath}
\begin{split}
\sum_n|e_n(z)\,e'_n(z)| &\leq
\left(\sum_n|e_n(z)|^2\right)^{1/2}\left(\sum_n|e'_n(z)|^2\right)^{1/2}
\\
& =
||K_z||_{\Ftwof}\left|\left|\frac{\partial}{\partial\bar{z}}K_z\right|\right|_{\Ftwof}.
\end{split}
\end{displaymath}
Now, the inequality \eqref{eqsp:kl} follows from Lemma
\ref{le:sp1}. This completes the proof of the Proposition.
\end{pf}

\par Now we turn to prove the necessity for  $0<p<\infty$.

\begin{proposition}\label{pr:Schnec}
Let $g\in H(\C)$, $0< p<\infty$ and $\phi\in\I$.
 If $T_g\in\mathcal{S}_p(\Ftwof)$, then $\frac{g'}{1+\phi'}\in L^p(\C, \Delta\phi\,dm)$.
\end{proposition}
\begin{pf}
We split the proof in two cases.
\par {\bf{ Case $\mathbf{2\le p<\infty}$.}}\,
 Suppose that $T_ g$ is in $\mathcal{S}_ p$, and let $\{e_ k\}$ be an
orthonormal basis in $\Ftwof$ and $R>100$. Let $\{z_ k\}$ be the sequence from Lemma
\ref{covering} for $t(z)=\delta\tau(z)$, $\delta\in (0,m_{\tau})$, and consider the operator $A$ taking $e_
k(z)$ to $f_{z_{k}}(z)=F_{z_ k,R}(z)/\tau(z_ k)$. It follows
from Proposition \ref{pr:fqp} that the operator $A$ is bounded on
$\Ftwof$.
Then $T_ g A$ belongs to
$\mathcal{S}_{p}$ (see \cite[p.$27$]{Zhu}), and by \cite[Theorem
$1.33$]{Zhu}
\begin{displaymath}
\sum_ k \|T_ g(f_{z_ k})\|_{\Ftwof}^p=\sum_ k \|T_ g A e_
k\|_{\Ftwof}^p<\infty.
\end{displaymath}
This together with  Proposition \ref{pr:BDK} and Theorem  \ref{Littlewood-Paley}
gives
\begin{equation*} \label{eq:nec1}
\begin{split}
\sum_ k & \frac{1}{\tau(z_
k)^p}
\left(\int_{D(\tau(z_k))}|g'(z)|^2\,(1+\phi'(z))^{-2}\,dm(z)\right
)^{p/2}
\\
&
\asymp \sum_ k \left(\int_{D(\tau(z_k))}|f_{z_
k}(z)|^2\,|g'(z)|^2\,(1+\phi'(z))^{-2}\,e^{-2\phi(|z|)} \,dm(z)\right )^{p/2}
\\
&
\lesssim \sum_ k\|T_ g(f_{z_ k})\|_{\Ftwof}^p
<\infty.
\end{split}
\end{equation*}
\par On the other hand, if $\delta$ is sufficiently small, applying Lemma \ref{le1},
Lemma \ref{tau}, Lemma \ref{almostct} and Lemma \ref{covering}, it follows that
\begin{displaymath}
\begin{split}
\int_{\C}|g'(z)|^p &\,(1+\phi'(z))^{-p}\,\Delta \phi(z) \,dm(z)
\\&\lesssim \sum_ k
\int_{D(\delta\tau(z_k))}\!\!\left(\frac{1}{\tau^2(z)}\int_{D(\delta\tau(z))}\!|g'(\zeta)|^2dm(\zeta)\right)^{p/2}
\!\!(1+\phi'(z))^{-p}\,\frac{dm(z)}{\tau^2(z)}
\\
&\lesssim \sum_ k\frac{1}{\tau (z_
k)^p}\int_{D(\delta\tau(z_k))}\!\!\!\left
(
\int_{D(\delta\tau(
z))}\!|g'(\zeta)|^2\,(1+\phi'(\zeta))^{-2}\,dm(\zeta)\right)^{p/2}\!\frac{dm(z)}{\tau^2(z)}
\\
&\lesssim \sum_ k \frac{1}{\tau(z_ k)^p}\left
(\int_{D(3\delta\tau(z_k))}\!
|g'(\zeta)|^2\,(1+\phi'(\zeta))^{-2}\,dm(\zeta)\right )^{p/2}.
\end{split}
\end{displaymath}
This together with
the previous inequality concludes the proof.
\par {\bf{ Case $\mathbf{0<p<2}$.}}\,
If $T_ g\in \mathcal{S}_ p$ then the positive operator $T_ g^*T_ g$
belongs to $\mathcal{S}_{p/2}$. Without loss of generality we may
assume that $g'\neq 0$. Suppose
$$T^*_ g T_ g f=\sum_ n \lambda_ n \langle f,e_ n \rangle_* \,e_ n$$
is the canonical decomposition of $T_ g^*T_ g$. Then a standard argument gives that $\{e_ n\}$
is  an orthonormal basis.
So relation (\ref{eqRKbasis})
 together with Corollary \ref{co:testfunc} and H\"{o}lder's
inequality yields
\begin{displaymath}
\begin{split}
\int_{\C}|g'(z)|^p \,(1+\phi'(z))^{-p}\,&\Delta \phi(z)\,dm(z)
\\
&\asymp
\int_{\C}|g'(z)|^p \,(1+\phi'(z))^{-p}\,\|K_
z\|^2\,e^{-2\phi(|z|)}\,dm(z)
\\&=
\sum_ n \int_{\C}|g'(z)|^p \,(1+\phi'(z))^{-p}\,|e_
n(z)|^2\,e^{-2\phi(|z|)}\,dm(z)
\\
&\leq \sum_ n \left (\int_{\C}|g'(z)|^2 \,(1+\phi'(z))^{-2}\,|e_
n(z)|^2\,e^{-2\phi(|z|)}\,dm(z)\right)^{p/2}
\\
&\lesssim \sum_ n \langle T_ g^*T_ g e_ n,e_ n \rangle_* ^{p/2} =\sum_
n \lambda_ n^{p/2}=\|T_ g^*T_ g\|_{\mathcal{S}_{p/2}}^{p/2}.
\end{split}
\end{displaymath}
This completes the proof.
\end{pf}
\par Finally, we shall prove the main result in this section.
\medskip\par\begin{Pf}{\emph{{\,\,Theorem \ref{th:Schatten}.}}}
\par Part $(a)$ follows directly from Proposition \ref{pr:Schsuf} and Proposition
\ref{pr:Schnec}. Moreover, if $0<p\le 1$ and $T_g\in
\mathcal{S}_p(\Ftwof)\subset \mathcal{S}_1(\Ftwof) $,  then Proposition \ref{pr:Schnec} and Lemma \ref{distortion} (d) imply that for some $r_0\in (0,\infty)$
\begin{displaymath}
\begin{split}
\int_{r_0<|z|<\infty}\frac{|g'(z)|}{|z|} \,dm(z)&\lesssim
\int_{r_0<|z|<\infty}\frac{|g'(z)|}{(1+\phi'(z))}\, \Delta\phi(z)\,dm(z)
\\ & \le \int_{\C}\frac{|g'(z)|}{|1+\phi'(z)|}\, \Delta\phi(z)\,dm(z)  <\infty.
\end{split}
\end{displaymath}
Therefore,  it follows that
$g'\equiv 0$, which gives $(b)$. The proof is complete.
\end{Pf}

\section{Examples}\label{examples}
In this section, several examples of rapidly increasing functions  are
given. We offer the corresponding description for the boundedness and
compactness of the
integration operator $T_ g$ in each case.

\textbf{Example 1:} The functions $\phi(|z|)=|z|^\a$, $\a>2$ belong to $\I$ with  $\Delta\phi\asymp|z|^{\a-2}$, so we obtain the following as a byproduct of Theorem \ref{th:Bq}.

\begin{corollary}\label{th:potential}
Let  $0<p,q<\infty$, $g$ an entire function, $\phi(|z|)=|z|^\a$, $\a>2$.
\begin{enumerate}
\item[$(I)$]
\begin{enumerate}
\item[$(a)$] Let $0<p\le q<\infty$. Then,
\begin{itemize}
\item  If $1+(\a-2)\left(1-\frac{1}{p}+ \frac{1}{q}\right)<0$, then  $T_ g:\Fpf\rightarrow \Fqf$ is bounded if and only if
  $g$ is constant.
 \item  If $1+(\a-2)\left(1-\frac{1}{p}+ \frac{1}{q}\right)\ge 0$,
 then  $T_ g:\Fpf\rightarrow \Fqf$ is bounded if and only if
   $g$ is a polynomial with  $$deg(g)\le 2+ (\a-2)\left(1-\frac{1}{p} +\frac{1}{q}\right).$$

                                                                                                                    \end{itemize}
\item[$(b)$] Let $0<p\le q<\infty$. Then,
\begin{itemize}
\item  If $1+(\a-2)\left(1-\frac{1}{p}+ \frac{1}{q}\right)\le 0$, then  $T_ g:\Fpf\rightarrow \Fqf$ is compact if and only if
  $g$ is constant.
 \item  If $1+(\a-2)\left(1-\frac{1}{p}+ \frac{1}{q}\right)> 0$, then  $T_ g:\Fpf\rightarrow \Fqf$ is compact if and only if
   $g$ is a polynomial with  $$deg(g)< (\a-1)-\left(\frac{1}{p}-\frac{1}{q}\right)(\a-2)+1.$$
\end{itemize}
   \end{enumerate}
    \item[$(II)$] Let $0<q<p<\infty$.
    \begin{itemize}
    \item If $q\le \frac{2p}{p(\a-1)+2}$, $T_ g:\Fpf\rightarrow \Fqf$ is bounded if and only if $g$ is constant.
    \item If $q> \frac{\frac{2p}{\alpha-1}}{p+\frac{2p}{\alpha-1}}$.
    The following conditions are
equivalent:
\begin{enumerate}
\item[$(a)$] $T_ g:\Fpf\rightarrow \Fqf$ is compact;

\item[$(b)$] $T_ g:\Fpf\rightarrow \Fqf$ is bounded;

\item[$(c)$] $g$ is a polynomial with $deg(g)<\a-\frac{2}{r}$, where $r=\frac{pq}{p-q}$.
\end{enumerate}
\end{itemize}
\end{enumerate}
\end{corollary}
 \par Moreover,  Theorem \ref{th:Schatten} shows that,
  $T_g\in S_p$ if and only if $p>\frac{\a}{\a-1}$ and
 $g$ is a polynomial with $deg (g)<\a(1-\frac{1}{p})$.

The above example illustrates that there does not exist a  value of $\beta=\frac{q-p}{pq}$ big enough such that \eqref{tgqbigp} implies that $g$ is constant
 for all $\phi\in\I$.

\textbf{Example 2:} The functions $\phi(|z|)=e^{\beta|z|}$, $\beta>0$ belong to $\I$ with  $\Delta\phi\asymp e^{\beta|z|}$, so we obtain the following as a byproduct of Theorem \ref{th:Bq}.

\begin{corollary}\label{th:potential}
Let  $0<p,q<\infty$, $g$ an entire function,  $\phi(|z|)=e^{\beta|z|}$, $\beta>0$.
\begin{enumerate}
\item[$(I)$]
\begin{enumerate}
\item[$(a)$] Let $0<p\le q<\infty$. Then,
\begin{itemize}
\item  If $\frac{1}{p}-\frac{1}{q}>1$, then  $T_ g:\Fpf\rightarrow \Fqf$ is bounded if and only if
  $g$ is constant.
 \item  If $\frac{1}{p}-\frac{1}{q}=1$, then  $T_ g:\Fpf\rightarrow \Fqf$ is bounded if and only if
   $g$ is a polynomial with  $deg(g)\le 1$.
   \item  If  $\frac{1}{p}-\frac{1}{q}<1$, then  $T_ g:\Fpf\rightarrow \Fqf$ is bounded if and only if
   $$sup_{z\in\C}|g'(z)|e^{\left(\frac{1}{p}-\frac{1}{q}-1\right)|z|}<\infty.$$                                                                                                                    \end{itemize}
\item[$(b)$] Let $0<p\le q<\infty$. Then,
\begin{itemize}
\item  If $\frac{1}{p}-\frac{1}{q}\ge 1$, then  $T_ g:\Fpf\rightarrow \Fqf$ is compact if and only if
  $g$ is constant.
   \item  If  $\frac{1}{p}-\frac{1}{q}<1$, then   $T_g$ is compact if the "little oh" versions of the boundedness conditions in $(a)$ hold.

                                                                                                      \end{itemize}
   \end{enumerate}
    \item[$(II)$] Let $0<q<p<\infty$.
    The following conditions are
equivalent:
\begin{enumerate}
\item[$(a)$] $T_ g:\Fpf\rightarrow \Fqf$ is compact;

\item[$(b)$] $T_ g:\Fpf\rightarrow \Fqf$ is bounded;

\item[$(c)$] $\int_\C \left(|g'(z)|e^{-\beta|z|}\right)^r\,dm(z)<\infty$, where $r=\frac{pq}{p-q}$.
\end{enumerate}

\end{enumerate}
\end{corollary}
\par Moreover,   Theorem \ref{th:Schatten} shows that $T_g\in S_p$ if and only if
$$g'(z)e^{-\b\left(1-\frac{1}{p}\right)|z|}\in L^p(\C,dm),\quad p>1.$$

The above example shows that for any $r>0$
 in Theorem \ref{th:Bq} $(II)$ (c)
 $g$ can grow exponentially.

\textbf{Example 3:} The function $\phi(|z|)=e^{e^{|z|}}$ belongs to $\I$ with  $\Delta\phi\asymp e^{2|z|+e^{|z|}}$, so
Theorem \ref{th:Bq} provides the following

\begin{corollary}\label{th:potential}
Suppose $g$ is an entire function,  $0<p,q<\infty$  and $\phi(|z|)=e^{e^{|z|}}$.
\begin{enumerate}
\item[$(I)$]
\begin{enumerate}
\item[$(a)$] Let $0<p\le q<\infty$. Then,
\begin{itemize}
\item  If $\frac{1}{p}-\frac{1}{q}\ge1$, then  $T_ g:\Fpf\rightarrow \Fqf$ is bounded if and only if
  $g$ is constant.
   \item  If $\frac{1}{p}-\frac{1}{q}<1$, then  $T_ g:\Fpf\rightarrow \Fqf$ is bounded if and only if
   $$sup_{z\in\C}|g'(z)|e^{|z|+\left(\frac{1}{p}-\frac{1}{q}-1\right)(2|z|+e^{|z|})}<\infty.$$                                                                                                                    \end{itemize}
\item[$(b)$] Let $0<p\le q<\infty$. Then $T_g$ is compact if the "little oh" versions of the boundedness conditions in $(a)$ hold.
\end{enumerate}
 \item[$(II)$] Let $0<q<p<\infty$.
    The following conditions are
equivalent:
\begin{enumerate}
\item[$(a)$] $T_ g:\Fpf\rightarrow \Fqf$ is compact;

\item[$(b)$] $T_ g:\Fpf\rightarrow \Fqf$ is bounded;

\item[$(c)$] $\int_\C \left(|g'(z)|e^{-|z|-e^{|z|}}\right)^r\,dm(z)<\infty$, where $r=\frac{pq}{p-q}$.
\end{enumerate}
\end{enumerate}
\end{corollary}
\par Moreover,   Theorem \ref{th:Schatten} shows that $T_g\in S_p$ if and only if$$g'(z)e^{\left(\frac{2}{p}-1\right)|z|-\left(1-\frac{1}{p}\right)\displaystyle{e^{|z|}}}\in L^p(\C,dm),\quad p>1.$$

\section{Invariant subspaces of the Volterra operator}\label{invariant}

In the particular case $g(z)=z$, the operator $T_g$ becomes the  Volterra operator
$$
Vf(z)=\int_0^z f(\z) \, d\z, \quad z\in\C,\quad f\in\Fpf,\  p>0.
$$
Recall that a closed subspace $\M$ in $\Fpf$
is called {\it invariant} for $V$ if $V\M\subseteq \M$.
In this section we aim to characterize  the invariant subspaces of $V$.
For this we begin by showing that polynomials are dense in $\Fpf$, by two different methods.
The first one is based on some smooth polynomials and the second one follows by a dilation argument.
\subsection{Density of polynomials}
For the first proof, we  need some background on certain smooth polynomials defined in terms of Hadamard products.
Let  $\T$ be the boundary of  the unit disc $\D=\{z:|z|<1\}$.
If $W(e^{i\theta})=\sum_{k\in J\subset \Z}b_k e^{ik\theta}$ is a trigonometric polynomial and
$f(e^{i\theta})=\sum_{k\in\Z}a_ke^{ik\theta}\in L^p(\T)$,  then the Hadamard
product
    $$
    (W\ast f)(e^{i\theta})=\sum_{k\in J}a_kb_ke^{ik\theta}
    $$
is well defined.

If $\Phi:\mathbb{R}\to\C$ is a $C^\infty$-function with compact
support $\supp(\Phi)$, we set
    $$
    A_{\Phi,m}=\max_{s\in\mathbb{R}}|\Phi(s)|+\max_{s\in\mathbb{R}}|\Phi^{(m)}(s)|,
    $$
and we consider the polynomials
    $$
    W_N^\Phi(e^{i\theta})=\sum_{k\in\mathbb
    Z}\Phi\left(\frac{k}{N}\right)e^{ik\theta},\quad N\in\N.
    $$
With this notation we can state the next result on smooth partial
sums.

\begin{lettertheorem}\label{th:cesaro}
Let $\Phi:\mathbb{R}\to\C$ be a $C^\infty$-function with compact
support $\supp(\Phi)$. Then the following
assertions hold:
\begin{itemize}
\item[\rm(i)] There exists a constant $C>0$ such that
    $$
    \left|W_N^\Phi(e^{i\theta})\right|\le C\min\left\{
    N\max_{s\in\mathbb{R}}|\Phi(s)|,
    N^{1-m}|\theta|^{-m}\max_{s\in\mathbb{R}}|\Phi^{(m)}(s)|
    \right\},
    $$
for all $m\in\N\cup\{0\}$, $N\in\N$ and $0<|\theta|<\pi$.
\item[\rm(ii)] If $0<p\le 1$ and $m\in\N$ with $mp>1$, there exists a constant $C(p)>0$ such that
    $$
    \left(\sup_{N}\left|(W_N^\Phi\ast f)(e^{i\theta})\right|\right)^p\le CA_{\Phi,m}
    M(|f|^p)(e^{i\theta})
    $$
for all $f\in H^p$. Here $M$ denotes the Hardy-Littlewood
maximal-operator
    $$
    M(|f|)(e^{i\theta})=\sup_{0<h<\pi}\frac{1}{2h}\int_{\theta-h}^{\theta+h}|f(e^{it})|\,dt.
    $$
\item[\rm(iii)] For each $p\in(0,\infty)$ and $m\in\N$ with $mp>1$, there exists a constant
$C=C(p)>0$ such that
    $$
    \|W_N^\Phi\ast f\|_{H^p}\le C A_{\Phi,m}\|f\|_{H^p}
    $$
for all $f\in H^p$.
\par Moreover, if $\Phi(0)=1$, then
\item[\rm(iv)] $\lim_{N\to\infty} \|f-W_N^\Phi\ast f\|_{H^p}=0$,\, for any $f\in H^p$, $0<p<\infty$.
 \item[\rm(iv)] $\lim_{N\to\infty}\left(W_N^\Phi\ast f\right)(e^{i\theta})=f(e^{i\theta})$ a.~e. $e^{i\theta}\in\T$, for any $f\in H^p$ and $0<p<\infty$.
\end{itemize}
\end{lettertheorem}

Theorem~\ref{th:cesaro} follows from the results and proofs in
\cite[p.~111-113]{Pabook}.

\par With this tool in our hands, we are going to prove the following.
\begin{theorem}\label{density}
Assume that
$\phi\in\I$ .
Then the polynomials are dense in $\Fpf$.
\end{theorem}
\noindent {\it Proof 1.}
Pick up a function $\Phi:\mathbb{R}\to\C$ be a $C^\infty$-function with $\Phi(0)=1$ and compact
support. For simplicity, we shall write $W_n$ for $W_n^\Phi$.
\par Take $f\in \Fpf$ and fix $r>0$. Let us denote $f_r(e^{i\theta})=f(re^{i\theta})$ and let us consider the sequence of polynomials
$\{u_n\}$, where
 $u_n(re^{i\theta})=\left(W_n\ast f_r\right)(e^{i\theta})$. By  $(iii)$ of Theorem \ref{th:cesaro}, there is a constant $C$ which does not depend on $r$ neither nor on $n$  such that
 $$M^p_p(r,u_n)\le C M^p_p(r,f)$$
 so for any $0<R<\infty$
 \begin{equation}\label{eq:dens1}
 \int_{|z|>R} |u_n(z)|^p e^{-p\phi(|z|)}\,dm(z)\le C \int_{|z|>R} |f|^p e^{-p\phi(|z|)}\,dm(z)\le ||f||^p_{\Fpf}<\infty.
 \end{equation}
 So $\sup_{n}||u_n||_{\Fpf}<\infty$, which together to Lemma \ref{le1} implies that $\{u_n\}$ is uniformly bounded on each compact subset of $\C$, so by Montel's theorem $\{u_n\}$ is a normal family. Consequently there is a subsequence  $\{u_{n_k}\}$ that converges to an entire function $g$ on compact subsets of $\C$. By $(iv)$ of Theorem \ref{th:cesaro} $g=f$. Finally, joining this fact to \eqref{eq:dens1}, and standard argument finishes the proof.
$\hfill\Box$

\medskip

\noindent {\it Proof 2.}
 For $f\in{\Fpf}$ and $r\in(0,1)$ put $f_r(z)=f(rz),\ z\in \C$.
Then
\begin{eqnarray}\label{dilationsf}
\lim_{r\rightarrow 1}\|f_r-f\|_{\Fpf}=0.
\end{eqnarray}
This follows by a standard argument that we sketch for the sake of completeness.
Let $R>0$ be such that $\int_{|z|>R} |f|^p e^{-p\phi} \, dA<\varepsilon$. Hence
\begin{eqnarray*}
\int_\C |f_r-f|^p e^{-p\phi} \, dm&\le& \int_{|z|\le R} |f_r-f|^p e^{-p\phi} \, dm + 2^p\int_{|z|>R} |f|^p e^{-p\phi} \, dm\\
&&\hspace{3.9cm}+2^p\int_{|z|>R} |f_r|^p e^{-p\phi} \, dm.
\end{eqnarray*}
It is clear that the first term in the above sum goes to zero as $r\rightarrow 1^-$.
Now using polar coordinates together the fact that the integral means $M_p(\rho,f)$ are increasing in $\rho$ we deduce that
$$
\int_{|z|>R} |f_r|^p e^{-p\phi}\, dm\le \int_{|z|>R} |f|^p e^{-p\phi}\, dm<\varepsilon,
$$
and now (\ref{dilationsf}) follows.

We now show that $f_r$ can be approximated by its Taylor polynomials in the ${\Fpf}-$norm.

Suppose first that $p\ge 1$.
If $a_n$ is the $n-$th coefficient in the Taylor expansion of $f$, then by the Cauchy formula and H\"older's inequality we get
$$
|a_n|^p r^{np}\le c_p \int_0^{2\pi} {|f(re^{it})|^p}\, dt,\quad r\ge 0.
$$
Multiplying both sides of the above inequality by $re^{-p\phi(r)}$ and integrating on $[0,\infty)$ we obtain
$$
|a_n| \|z^n\|_{\Fpf} \le c_p \|f\|_{\Fpf}.
$$
Using this we obtain
$$
\left\|\sum_{n=N}^\infty a_n r^n z^n \right\|_{\Fpf}\le \sum_{n=N}^\infty |a_n| \|z^n\|_{\Fpf} r^n\le c_p \|f\|_{\Fpf} \sum_{n=N}^\infty r^n,
$$
which shows that $f_r$ can be approximated by its Taylor polynomials in the ${\Fpf}-$norm.

Assume now $0<p<1$.   We have  (see \cite[Theorem $6.4$]{Duren1970})
$$
 r^{np}|a_n|^p \le C_p \, n^{1-p}\int_0^{2\pi} |f(re^{it})|^p dt,
$$
and proceeding as in the case $p\ge 1$ we deduce
$$
|a_n| \|z^n\|_{\Fpf} \le C_p n^{1/p-1}\|f\|_{\Fpf}.
$$
Thus
$$
\|\sum_{n=N}^\infty a_n r^n z^n \|^p_{\Fpf}\le \sum_{n=N}^\infty |a_n|^p \|z^n\|^p_{\Fpf} r^{np}\le C_p\|f\|_{\Fpf} \sum_{n=N}^\infty r^{np} n^{1-p},
$$
and the conclusion follows, since the last expression above is the tail of a convergent series.
$\hfill\Box$

\subsection{Invariant subspaces}
Notice that, with respect to the standard orthonormal basis $e_n(z)=z^n/\|z^n\|_{\Ftwof}$, $n\ge 0$, the operator $V:\Ftwof\rightarrow\Ftwof$ is a weighted shift, i.e. $Ve_n=\om_n e_{n+1}$ with
weight sequence
\begin{eqnarray}\label{weiseq}
\om_n=\frac{\|z^{n+1}\|_{\Ftwof}}{(n+1)\|z^n\|_{\Ftwof}},\quad n\ge 0.
\end{eqnarray}
A simple argument using integration by parts and Lemma \ref{distortion} shows that, for $\phi\in \I$, we have $\om_n\rightarrow 0$ as $n\rightarrow\infty$.

\begin{theorem}\label{invs}
Assume $p>0$ and let $\phi\in\I$ 
 If the sequence $(\om_n)_{n\ge 1}$ given by (\ref{weiseq}) is eventually decreasing to zero, then a closed subspace $\M\subset \Fpf$ is a proper invariant subspace of $V:\Fpf\rightarrow \Fpf$
if and only if  there exists a positive integer $N$ such that
$$
\M=\{f\in\Fpf\,:\, f^{(k)}(0)=0 \hbox{ for } 0\le k\le N-1\}=\overline{\hbox{Span}\{z^k\,:\, k\ge N \}}^{\Fpf}.
$$
\end{theorem}

\begin{pf}
Clearly, the sets
$$
A_N^p:=\{ f \in\Fpf\,:\, f^{(k)}(0)=0 \hbox{ for } 0\le k\le N-1\}
$$
are invariant subspaces for $V$.

Let us now prove that these are all the invariant subspaces of $V$.
The case $p=2$ follows directly from a result of Yakubovich on weighted shifts (see \cite{yakubovic} Theorems 3-4).
We are now going to show that the result for $p\neq 2$ can be obtained via the case $p=2$ using the boundedness of $T_g$ on Fock spaces
$\Fpf$ (see Theorem \ref{th:Bq} above).
In order to do this, let us first prove the following

{\it Claim:} For each $p>0$ there exist positive integers $M_1=M_1(p),M_2=M_2(p)$ such that $V^{M_1}:\Ftwof\rightarrow\Fpf$
and $V^{M_2}:\Fpf\rightarrow\Ftwof$ are bounded.

Let us first notice that by Theorem \ref{th:Bq} and by $(a), (b)$ of Lemma \ref{distortion}  we have

\begin{eqnarray}
V:\mathcal{F}_{p'}^\phi \rightarrow \mathcal{F}_{q}^\phi \hbox{ is bounded if } F(q)\le p'\le q\label{inv},\\
V:\mathcal{F}_{p'}^\phi \rightarrow \mathcal{F}_{q}^\phi \hbox{ is bounded if } F(p')<q< p'\label{inv0},
\end{eqnarray}
where $F(s)=\frac{2s}{s+2}$ for $s>0$. From this it is easy to see that $V:\Fpf\rightarrow\Ftwof$ and $V:\Ftwof\rightarrow \Fpf$ are bounded
if $p\ge 2$, and hence we can take $M_1=M_2=1$ in this case.

Let us now study the case $0<p<2$. Consider  sequence $(x_n)$ given by $x_0=2$, $x_{n}=F(x_{n-1})$, $n\ge 1$ and notice that
$(x_n)$ decreases to zero as $n\rightarrow\infty$.
We aim to show that
 \begin{eqnarray}
 &V^M:\Fpf\rightarrow \Ftwof\hbox{ is bounded if }& x_{M-1}\ge p \ge x_M,\label{inv1}\\
 &V^M:\Ftwof \rightarrow \Fpf\hbox{ is bounded if }& x_{M-1}\ge p > x_M,\, M\ge 1\label{inv2}.
 \end{eqnarray}
 Both assertions follow by induction. We only prove (\ref{inv2}), since (\ref{inv1}) is shown in a similar way.
For $M=1$ relation ($\ref{inv2}$) follows directly from (\ref{inv}-\ref{inv0}).
Now assume that $V^M:\Ftwof\rightarrow\Fpf$ is bounded if  $x_{M-1}\ge p > x_M$. By (\ref{inv0}) we deduce that
$V: \mathcal{F}_{x_{M}+\varepsilon}^\phi\rightarrow\Fqf$ is bounded if
$
x_{M}+\varepsilon>q>F(x_{M}+\varepsilon).
$
For $0<\varepsilon<x_{M-1}-x_M$ we now obtain
$$
V^{M+1}:\Ftwof\rightarrow\Fqf \hbox{ is bounded if } x_{M}+\varepsilon>q>F(x_{M}+\varepsilon).
$$
Let $\varepsilon\rightarrow 0$ above to get
$$
V^{M+1}:\Ftwof\rightarrow\Fqf \hbox{ is bounded if } x_{M}\geq q> F(x_{M})=x_{M+1},
$$
and with this the claim is proven.

Now let $\M$ be an invariant subspace for $V:{\mathcal{F}^\phi_p}\rightarrow{\mathcal{F}^\phi_p}$.
Then let $M_1,M_2\in\N$ be such that $V^{M_1}:\Ftwof\rightarrow\Fpf$ and $V^{M_2}:\Fpf\rightarrow\Ftwof$
are bounded.
Then $V^{M_2}\M \subset {\mathcal{F}^\phi_2}$ and $\overline{V^{M_2}\M}^{\mathcal{F}^\phi_2}$ is an invariant subspace for $V$ on ${\mathcal{F}^\phi_2}$,
so it is of the form
$$\overline{V^{M_2}\M}^{\mathcal{F}^\phi_2}=A^2_N.$$ Now let $f\in A^2_N$. Then there exists a sequence $(f_n)\subset \M$ with
$$
V^{M_2}f_n\rightarrow f \ \ in \ \ \mathcal{F}^\phi_2.
$$
But $V^{M_1}:\mathcal{F}^\phi_2\rightarrow \mathcal{F}^\phi_p$ is bounded, hence
$$
V^{{M_1+M_2}} f_n\rightarrow V^{M_1}f \ \ in \ \ \mathcal{F}^\phi_p,
$$
which gives $V^{M_1} f\in \M$, since $V^{{M_1+M_2}} f_n\in \M$. Put $f=z^N$ to deduce $z^{N+M_1}\in \M$ and by applying $V$ indefinitely to $z^{N+M_1}$  we get
$$
A^p_{N+M_1}\subseteq \M.
$$
If $\M \neq A^p_{N+M_1}$, let $N_1$ be the smallest nonnegative integer such that there exists $f\in\M$ with
$f^{(N_1)}(0)\neq 0$ (clearly $0\le N_1<N+M_1$). But, for this particular $f$, we have $V^{N+M_1-N_1-1}f \in \M$, which
in view of the above inclusion implies $z^{N+M_1 -1}\in\M$, and therefore $A^p_{N+M_1-1}\subseteq \M$. We repeat this procedure until we obtain
$A^p_{N_1}\subseteq \M$, and then the choice of $N_1$ forces $\M=A^p_{N_1}$, so that the proof
is done.
\end{pf}

\begin{corollary}
For $\phi(z)=|z|^\alpha$ with $\alpha>2$, the proper invariant subspaces of $V$ on $\Fpf$ are precisely the spaces
$$
\{f\in\Fpf\,:\, f^{(k)}(0)=0 \hbox{ for } 0\le k\le N-1\}=\overline{\hbox{Span}\{z^k\,:\, k\ge N \}}^{\Fpf}.
$$
\end{corollary}

\begin{pf}
By Theorem \ref{invs} it suffices to show that  the sequence
$$
\om_n=\frac{\|z^{n+1}\|_{\Ftwof}}{(n+1)\|z^n\|_{\Ftwof}},\quad n\ge 0
$$
is eventually decreasing to zero. By Stirling's formula we deduce that
$$
\om_n^2= 2^{-2/\alpha}\frac{\Gamma(\frac{2}{\alpha}(n+2))}{(n+1)^2\Gamma(\frac{2}{\alpha}(n+1))} \asymp n^{2/\alpha -2},
$$
and hence $\om_n\rightarrow 0$ as $n\rightarrow\infty$. From the proof of Proposition 7 in \cite{atzmon} it follows that
the function
$$
f(x)= \frac{\Gamma(x+2/\alpha)}{x^2 \Gamma(x)}
$$
is eventually decreasing for $\alpha>1$. Thus $\{\om_n\}_{n=0}^\infty$ is eventually decreasing to zero.
\end{pf}
\section{Remarks on the Bergman space case}\label{bergman}
In this section, we would like to illustrate that
some of the methods we employed in the context of Fock spaces
provide additional insight into the above mentioned results from \cite{PP,PavP}.

 Given a positive radial weight $w$,
 the distortion function is defined  as follows (see \cite{PP,PavP})
$$
\psi_w(r)=\frac{1}{w(r)} \int_r^1 w(u)\, du,\quad 0\le r<1.
$$
Analogues of Theorem \ref{th:Bq} and Corollary \ref{co:LPf} for the setting of Bergman spaces with rapidly decreasing weights were
obtained in \cite{PP,PavP}.  In particular, the following holds.
\begin{lettertheorem}\label{PavP}
Suppose that  $w$ is  a radial differentiable weight, and there is $L>0$
such that
\begin{eqnarray}\label{eq:L}
\sup_{0<r<1}\frac{w'(r)}{w(r)^2}\int_r^1 w(x)\,dx\ \le L,
\end{eqnarray}
 then for each $p\in (0,\infty)$ and $g$ analytic on $\D$
\begin{displaymath}\label{LPformulaBergman}
\int_{\D}| g(z)|^p w(z)\,dm(z)\asymp
|g(0)|^p+\int_{\D}|g'(z)|^p\,\psi_ w(z)^p\,w(z)\,dm(z).
\end{displaymath}
\end{lettertheorem}

Let
$$
\phi(r)=-\log w(r), \quad 0\le r<1.
$$
The following disc analogue of Lemmas \ref{lhospital}, \ref{distortion} and \ref{almostct} can be easily proven.
\begin{lemma}\label{lhospitalbergman}
Assume $\phi:[0,1)\rightarrow \R$ is twice continuously differentiable and  there exists $r_0\in[0,1)$ such that $\phi'(r)\neq 0$ for $1>r>r_0$.
Suppose
\begin{eqnarray}\label{eqnb}
&&\lim_{r\rightarrow 1} \frac{e^{-\phi(r)}}{\phi'(r)}=0 \nonumber\\
&&\liminf_{r\rightarrow 1} \frac{\phi''(r)}{\phi'(r)^2} >-1\nonumber\\
&&\limsup_{r\rightarrow 1}\frac{\phi''(r)}{\phi'(r)^2}<\infty.\nonumber
\end{eqnarray}
Then \eqref{eq:L} holds and there exists $r_1\in [0,1)$ such that
\begin{eqnarray}\label{brg}
\psi_w(r)\asymp \frac{1}{\phi'(r)} \hbox{ for } r\in [r_1,1)
\end{eqnarray}

\end{lemma}
\begin{pf}
By hypotheses there is $\alpha>-1$ and  $r_2\ge r_0$ such that $\frac{\phi''(r)}{\phi'(r)^2}\ge \a$ on $[r_2,1)$. So,
an integration by parts  on $(r_2,r]\subset (r_2,1)$, gives
\begin{eqnarray*} \label{eq:ip}
&\int_{r_2}^r e^{-\phi(s)}\,ds= \int_{r_2}^r \frac{-\phi'(s)e^{-\phi(s)}}{-\phi'(s)}\,ds
 =\frac{-e^{-\phi(r)}}{\phi'(r)}+\frac{e^{-\phi(r_2)}}{\phi'(r_2)}-\int_{r_2}^r \frac{\phi''(s)}{\phi'(s)^2} e^{-\phi(s)}\,ds
\\ & \le \frac{-e^{-\phi(r)}}{\phi'(r)}+\frac{e^{-\phi(r_2)}}{\phi'(r_2)}-\alpha\int_{r_2}^r  e^{-\phi(s)}\,ds,
\end{eqnarray*}
that is
\begin{eqnarray*}
&\int_{r_2}^r e^{-\phi(s)}\,ds\le\frac{1}{\alpha+1}\left(\frac{-e^{-\phi(r)}}{\phi'(r)}+\frac{e^{-\phi(r_2)}}{\phi'(r_2)}\right)
\end{eqnarray*}
so taking limits as $r\to 1^-$, we deduce that $\int_0^1 e^{-\phi(s)}\,ds<\infty$. Next, arguing as in the proof Lemma \ref{lhospital}, \eqref{brg} follows.
These calculations also give \eqref{eq:L}. This finishes the proof.
\end{pf}

\par It is worth noticing that there are weights $\om$ satisfying \eqref{eq:L} but such that \eqref{brg} does not hold.
 The weight
 $\om(r)=e^{-(1-r)^\a}\asymp 1$,  $\a>0$,  gives a concrete example.
\begin{lemma}\label{distortionbergman}
Assume  $\phi:[0,\infty)\rightarrow \R^+$ is a twice continuously differentiable function such that $\Delta\phi>0\, ,(\Delta \phi (z))^{-1/2} \asymp \tau (z)$, where $\tau(z)$ is a radial positive function that
decreases to zero as $|z|\rightarrow 1^-$ and $\lim_{r\rightarrow 1^-} \tau'(r)=0$.
Then
\begin{eqnarray*}
&(a)& \lim_{r\rightarrow 1^-} (1-r)\phi'(r)=\infty.\\
&(b)& \lim_{r\rightarrow 1^-}  \tau(r) \phi'(r)=\infty, \hbox{or, equivalently, } \lim_{r\rightarrow 1^-} \frac{\phi''(r)}{(\phi'(r))^2}=0. \\
&(c)&
\psi_\om(r)\asymp \frac{1}{\phi'(r)+1}, \hbox{ for }
r \in[ 0,1).\\
&(d)& \hbox{There exists $r_0\in[0,1)$ such that  for all $a\in\D$ with $1>|a|>r_0$, and  any}
\\  &&\hbox{$\delta>0$ small enough we have}\end{eqnarray*}
$$\phi'(a)\asymp\phi'(z), \quad z\in D(a, \delta \tau(a)).$$

\end{lemma}

\begin{pf}
Since by L'Hospital's rule
$$\lim_{r\rightarrow 1^-} \frac{1-r}{\tau(r)}
=\lim_{r\rightarrow 1^-} \frac{-1}{\tau'(r)}=+\infty$$
 bearing in mind (\ref{laplace})  and using  again L'Hospital's rule,
\begin{equation*} \begin{split}
&\lim_{r\rightarrow 1^-}r(1-r)\phi'(r)=\lim_{r\rightarrow 1^-}\frac{\int_0^r  s\Delta \phi (s)\,ds}{(1-r)^{-1}}
\\ &\ge C\lim_{r\rightarrow 1^-}\frac{\int_0^r \frac{s}{\tau^2(s)}\, ds}{(1-r)^{-1}}
\ge C\lim_{r\rightarrow 1^-} \frac{r(1-r)^2}{\tau^2(r)}=
+\infty,
\end{split}\end{equation*}
which gives $(a)$.

Arguing as in the proof of Lemma \ref{distortion}
we obtain
$\lim_{r\rightarrow 1^-}\tau(r)\phi'(r)=\infty$. By $(a)$ and relation (\ref{laplace}) this last
fact is equivalent to $\lim_{r\rightarrow 1^-}\frac{ \phi''(r)}{(\phi'(r))^2}=0$.

Taking into account $(a)-(b)$ it is clear that the hypotheses in Lemma \ref{lhospitalbergman} are satisfied,
 and hence
$\psi_\om(r)\asymp \frac{1}{\phi'(r)}$ for $r\ge [r_0,1)$. Since $\phi'\ge 0$ and $\lim_{r\rightarrow 1^- } \phi'(r)=\infty$, we obtain
$(c)$.
\par Part (d) can be proved following the  steps  in the proof of Lemma \ref{almostct}.
\end{pf}

As a byproduct of the previous lemmas, we obtain the following improvements:
\begin{itemize}
\item In view of (\ref{brg}) and \cite[Theorem B]{PP}
the more transparent Littlewood-Paley formula
\begin{eqnarray}\label{LPrdberg}
||f||^p_{A^p_\om}\asymp |f(0)|^p+\int_\D \frac{|f'(z)|^p}{(1+\phi'(z))^p}\,\om(z)dA(z)
\end{eqnarray}
 can be obtained for Bergman spaces with weights $\om=e^{-\phi}$
belonging to the class $\I$ considered in \cite{PP}.

\item
Part (d) of Lemma \ref{distortionbergman}
shows that the hypothesis $(6)$ on the distortion function in  \cite[Theorem $2$]{PP} can be omitted.
 Going further, this observation allows us to
 extend the description of the boundedness and compactness of $T_g$  to a wider class of weights. Especially, we admit for a considerably faster decay, including
triple exponential weights of the form
$$
\omega(z)=\exp({-e^{e^\frac{1}{1-|z|}}}).
$$
\end{itemize}

\end{document}